# The structural equations of Cartan and the wave solution Einstein's equation

Vladimir Shcherban' (Scientific leader - Olga Baburova)

Moscow Pedagogical State University, Moscow, Russia

e-mail: shchvova@orc.ru

Work consists of introduction, two chapters, conclusion and four applications. In this work is examined the condition, with which the wave space metrics of Riemann- Cartan is the solution of Einstein's equation in the void. Geometric structures were for this purpose studied on the differentiated variety: connectedness, curvature and twisting connectedness. For this was used this analytical apparatus of differential geometry as the calculation of exterior forms.

In the first chapter the following concepts are examined:

- variety;

- vectors and 1- form's on the varieties;

- metric tensor, connectedness and the covariant derivative;

- form with the values in the vector spaces;

- the structural equations of Cartan.

The second chapter is dedicated to the presence of condition, during which wave certificate is the solution of Einstein's equation in the void. Using the first and second structural equations of Cartan for the variety without twisting of that allotted by wave certificate, they were calculated: the 1- form's of connectedness and the coefficients of connectedness, and also Riemann's tensor, Ricci's tensor and scalar of curvature.

In appendix 1 the derivation of formula for the coefficients of connectedness in the space with the twisting is given.

In appendix 2 the calculation of the components of the 1- forms of connectedness is carried out.

In appendix 3 the calculations of the components of the 2- forms of curvature are presented.

In appendix 4 the components of Ricci's tensor are calculated.

**Оглавление**

# Введение

Геометрия как наука, возникла более чем две тысячи лет назад и обязана своим существованием таким первичным феноменам мира, как пространство и время. Сегодня геометрия – это прежде всего дифференциальная геометрия, которая представляет собой один из самых важных разделов современной математики.

Первоначально созданная на основании координатного метода Р.Декарта и П.Ферма с учетом идей Н.И. Лобачевского, дифференциальная геометрия начала бурно развиваться благодаря вкладу Б.Римана, который сделал принципиально новый шаг – стал изучать произвольные, так называемые римановы пространства, сразу же нашедшие важные приложения в механике и в теории относительности Эйнштейна, которая в то время была отнюдь не бесспорна. Возникли новые разделы геометрии – векторное и тензорное исчисление, риманова геометрия, закладываются основы топологии, как части геометрии, посвященной свойству непрерывности. Большую роль в создании топологии и ее применении сыграли Л. Эйлер, К. Жордан, М. Френе, Ф. Хаусдорф, А. Лебег, Л. Бауер, А. Пуанкаре. Соединение идей нелинейных координат, векторного и тензорного исчисления, теории поверхностей, а также плодотворность геометрических представлений в различных задачах естествознания привели к важнейшему в геометрии понятию – многообразию.

Предметом дифференциальной геометрии является изучение геометрических задач, допускающих использование анализа. Поэтому объектами изучения являются пространства, в которых имеют смысл такие понятия как дифференцирование и интегрирование. Такой класс пространств образуют дифференцируемые многообразия. Важнейшим аналитическим методом изучения этих пространств является исследование их бесконечно малых участков. Основное преимущество перехода к бесконечно малым объектам заключается в том, что при этом все они становятся линейными. Каждая кривая в бесконечно малом есть прямая линия, в том смысле, что ее можно заменить касательной. Одним из самых полезных и плодотворных аналитических ме-тодов в дифференциальной геометрии является исчисление дифференциальных функций, созданный в начале XX века Э. Картаном.

Дополнительными структурами на многообразии являются связность и метрика, которые позволяют определить правило параллельного переноса и понятие расстояния, имеющие огромное прикладное значение в физике, в частности, в физике элементарных частиц.

Современная общая теория относительности, в частности – теория гравитации Эйнштейна, основана на геометризации физических законов. В качестве объекта, на котором строится геометрическая интерпретация теории, рассматривается многообразие пространства-времени – четырехмерное псевдориманово многообразие с лоренцевой метрикой.

В общей теории относительности гравитационное поле отождествляется с геометрической структурой пространственно-временного континуума на основе представления римановой геометрии. Известно, что все материальные тела прочерчивают в пространстве-времени одинаковые траектории, если на них действуют лишь гравитационные силы. На основании этого факта была выдвинута идея о том, гравитационное поле является неотъемлемым свойством пространства. Современные физики-релятивисты трактуют гравитационное поле, как геометрическую структуру пространства-времени – кривизну. Одной из основных задач современной теории гравитации является проблема гравитационных волн. Гравитационная волна – это возмущение гравитационного поля, то есть кривизны пространства-времени, распространяющееся с конечной скоростью и несущее с собой энергию.

В связи с выше сказанным, первая глава данной работы посвящена рассмотрению основных объектов и операций внешнего дифференциального исчисления Картана, обобщенного на случай произвольного дифференцируемого многообразия. Для этого необходимо изложить основы теории дифференцируемых многообразий. Кривые и функции на многообразии вводятся на языке координатных отображений, векторы трактуются как операторы дифференцирования вдоль кривых. Также определяются 1-формы, тензоры и p-формы, внешнее произведение форм, вводятся основные геометрические структуры на многообразии: метрика, связность, кривизна и кручение связности. Здесь же излагается на языке внешних форм основные геометрические уравнения на многообразии (в котором заключены все геометрические свойства последнего) - структурные уравнения Картана для кривизны.

Вторая глава работы посвящена изучению проблемы гравитационных волн. Известно, что гравитационные волны вдали от излучающих источников ведут себя как плоские электромагнитные волны, поэтому будем рассматривать плоские волны с заданной метрикой.

Работая над созданием общей теории относительности, Эйнштейну удалось установить связь между гравитацией и римановым пространством-временем. А значит, появилась возможность отождествить геометрические структуры с элементами гравитационного поля, что и сделал Эйнштейн, получив следующее уравнение:

$$G_{\mu\nu} = R_{\mu\nu} - \frac{1}{2}g_{\mu\nu}R = \chi T_{\mu\nu},$$

получившее название уравнения Эйнштейна.

В данной работе будем рассматривать так называемое пустое пространство (это пространство, где кроме гравитационного нет других тел и полей). В таком пространстве уравнение Эйнштейна заметно упрощается и приобретает вид

$$R_{\mu\nu} = 0.$$

Наша задача будет состоять в нахождении условия, при котором волновая метрика риманова пространства является решением уравнения Эйнштейна в пустоте. Для этой цели будут вычислены 1-формы связности, 2-формы кривизны и тензор Риччи, приравнивая который к нулю будет найдено искомое условие.

# Глава I

## Структурные уравнения Картана

### §1.1 Понятие многообразия

Обозначим через $R^n$ совокупность всех наборов из $n$ вещественных чисел $(x^1, x^2, \ldots, x^n)$.

Говорят, что некоторое множество $M$ наделено топологической структурой, если каждому его элементу, называемому точкой, соответствует семейство множеств, называемых окрестностями этого элемента, и удовлетворяющих следующей системе аксиом:

1. $x \in M$ содержится в каждой своей окрестности, $x \in U(x) \quad U(x) \subset M$.
2. Для всякой точки $y \in U(x)$ существует $U(y) \subset U(x)$.
3. Если заданы окрестности $U_1(x)$ и $U_2(x)$ то существует $U_3(x) = U_1(x) \cap U_2(x)$.
4. (аксиома отдельности или аксиома Хаусдорфа)

    Если точка $x$ не совпадает с точкой $y$, то существуют окрестности $U(x)$ и $U(y)$: $U(x) \cap U(y) = \emptyset$.

Множество $M$, наделенное топологической структурой, называется топологическим пространством.

***Многообразие*** – это топологическое пространство, каждая точка которого обладает окрестностью $U$, гомеоморфной области в $R^n$.

Поставим в соответствие каждой точке $P \in M$ набор чисел $(x^1, x^2, \ldots, x^n) \in R^n$, то есть вводим координаты для точки.

$$f: P \to \{x^\alpha\}, \quad \alpha = 1, 2, \ldots, n$$
$$f(P) = (x^1, x^2, \ldots, x^n)$$

Отображение $f$ называется *координатным*, $f$ - непрерывно и взаимно-однозначно.

Рассмотрим две пересекающиеся окрестности точки $P \in M$ (см. **Рис.1**)

**Рис1**.

Окрестности $U$ и $V$ из $M$ пересекаются. Соответствующие координатные отображения $f$ и g в $R^n$ дают в зоне пересечения два различных отображения (а, следовательно, две системы координат). Соответствие между этими координатами характеризует класс гладкости многообразия.

$$f: P \to f(P) = \{x^\alpha\}, \quad \alpha = 1, 2, \ldots, n$$
$$g: P \to g(P) = \{y^\beta\}, \quad \beta = 1, 2, \ldots, n$$

$f$ - координатное отображение, следовательно, существует $f^{-1}$.

Рассмотрим $g \circ f^{-1}: R^n \to R^n$, значит, имеем

$$\begin{cases} y^1 = y^1(x^1, x^2, \ldots, x^n) \\ y^2 = y^2(x^1, x^2, \ldots, x^n) \\ \ldots\ldots\ldots\ldots\ldots\ldots\ldots\ldots \\ y^n = y^n(x^1, x^2, \ldots, x^n) \end{cases}, \qquad (1.1.1а)$$

то есть

$$y^\beta = y^\beta(x^\alpha). \qquad (1.1.1в)$$

Преобразования (1.1.1) называются *функциями перехода*. Если эти функции $k$ раз дифференцируемы, то они называются функциями класса $C^k$.

Если $k = 1$, то есть существуют первые производные, то функции перехода (1.1.1) непрерывны.

**Определение**. Пара $(U, f)$ называется *картой*, где $U$ – окрестность, $f$ - координатное отображение.

Если функции перехода между картами многообразия $M$ являются функциями класса $C^k$, то такие карты называются $C^k$-согласованными.

**Определение.** Объединение всех карт многообразия $M$ называется *атласом*.

Топологическое пространство с его $C^k$-согласованными атласами называется многообразием (дифференцируемым многообразием при $k > 1$).

**Кривые на многообразиях.**

*Кривая* – это дифференцируемое отображение открытого подмножества из $\mathbb{R}^1$ в $M$ (см. рис. 2). Таким образом, каждой точке из $\mathbb{R}^1$ (являющимся вещественным числом, которое мы обозначим через $\lambda$) отвечает точка из $M$, называемая ***точкой-образом***. Вещественное число $\lambda$ называется ***параметром кривой***. Совокупность точек-образов и есть кривая. Если две кривые имеют совпавшие образы, а одинаковые точки отвечают различным значениям параметров, то эти кривые являются различными.

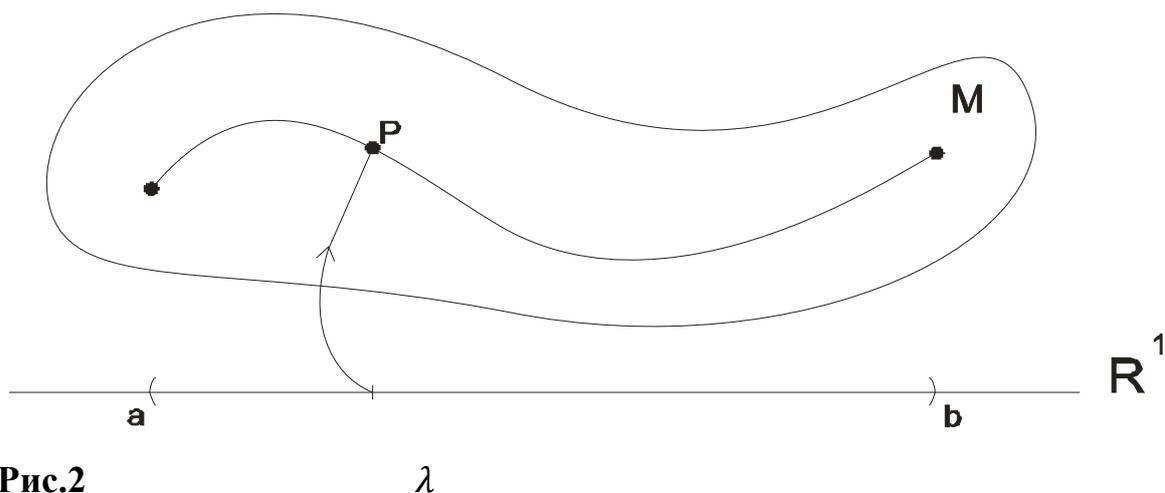

**Рис.2**

(Кривая в $M$ - это отображение из $\mathbb{R}^1$ в $M$. Точка $\lambda$ из $\mathbb{R}^1$ переходит в точку $P$ из $M$.)

**Функция на многообразиях.**

***Функция*** на *M* - это правило, сопоставляющее каждой точке из *M* вещественное число, которое называется ***значением функции***. Если некоторая область в *M* отображается при помощи гладкого координатного отображения на область в $\mathbb{R}^n$, то функция на *M* переходит в функцию на $\mathbb{R}^n$ (см. **рис. 3**). И если функция на $\mathbb{R}^n$ является дифференцируемой, то она ***дифференцируема*** на *M* [3].

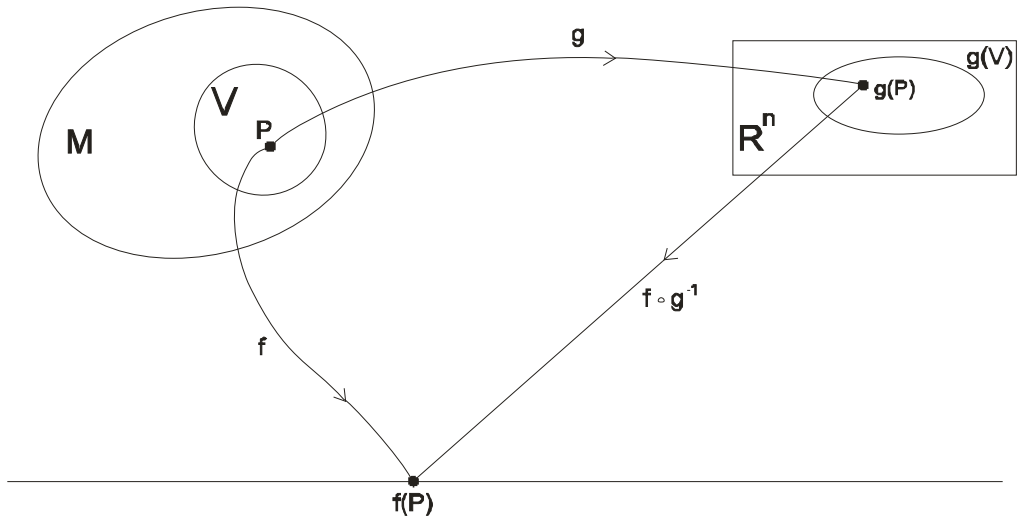

**Рис.4**. (Функция $f$ на *M* есть отображение из *M* в $\mathbb{R}^1$. Координатное отображение $g$ из области *V* многообразия *M*, содержащей *P*, на область $g(V)$ в $\mathbb{R}^n$ обратимо. Составное отображение $f \circ g^{-1}$ переводит $\mathbb{R}^n$ в $\mathbb{R}^1$, т.е. является функцией на $\mathbb{R}^n$. Это просто выражение $f(P)$ через координаты точки *P*).

### §1.2 Векторы и 1-формы на многообразиях

Определим понятие вектора на многообразии. Для этого рассмотрим кривую, которая задается уравнением вида $x^\alpha = x^\alpha(\lambda)$, $\alpha = 1, \ldots, n$ и дифференцируемую функцию $f(x^1, x^2, \ldots, x^n)$ в окрестности некоторой точки $P \in M$. В каждой точке кривой определено значение $f$. При этом на кривой возникает новая функция $g(\lambda)$ и

$$g(\lambda) = f(x^1(\lambda), x^2(\lambda), \ldots, x^n(\lambda)) = f(x^\alpha(\lambda)).$$

Дифференцируя ее получаем:

$$\frac{dg}{d\lambda} = \sum_\alpha \frac{dx^\alpha}{d\lambda} \frac{\partial f}{\partial x^\alpha} \qquad (1.2.1)$$

Это равенство справедливо для любой функции $g$, значит можно записать

$$\frac{d}{d\lambda} = \sum_\alpha \frac{dx^\alpha}{d\lambda} \frac{\partial}{\partial x^\alpha} \qquad (1.2.2)$$

Набор чисел $\left\{\frac{dx^\alpha}{d\lambda}\right\}$ в евклидовом пространстве является компонентами некоторого касательного вектора к кривой $x^\alpha(\lambda)$. Это следует из того, что $\{dx^\alpha\}$ – бесконечно малое перемещение вдоль кривой и деление его на $d\lambda$ меняет лишь длину, а не направление этого перемещения. Поскольку кривая берется с определенным параметром, для каждой кривой однозначно определен набор $\left\{\frac{dx^\alpha}{d\lambda}\right\}$. Значит, каждая кривая имеет единственный касательный вектор. Поэтому в точке $P$ существует касательный вектор к этой кривой, а так как в точке $P$ можно провести бесконечное множество кривых заданного типа, то возможно и бесконечное множество касательных векторов к этим кривым. Множество касательных векторов образует пространство $T_p$, которое называется **касательным пространством ковариантных векторов** к $M$ в точке $P$.

Из вышесказанного следует, что в теории дифференцируемых многообразий касательный вектор в точке $P \in M$ является линейным дифференциальным оператором, действующим на множество гладких функций $f$, заданных на $M$ в окрестности точки $P$. Операции дифференцирования вдоль координатных линий имеют вид частных производных $\frac{\partial}{\partial x^\alpha}$. Система $\left\{\frac{\partial}{\partial x^\alpha}\right\}$ образует **координатный базис** векторного пространства. Пространство всех касательных векторов $T_p$ в точке $P \in M$ и пространство всех дифференцируемых функций вдоль кривых, проходящих через точку $P$, находятся во взаимно-однозначном соответствии. Касательные векторы в точке $P$ можно выразить не только через координатный базис, но и через произвольный базис $\{\bar{e}_\alpha\}$

$$\bar{v} = \sum_\alpha v^\alpha \frac{\partial}{\partial x^\alpha} \qquad \bar{v} = v^\alpha \bar{e}_\alpha \qquad (1.2.3)$$

Правило, задающее вектор в каждой точке многообразия, называется **векторным полем**. Поэтому в точке $P$ касательные вектора образуют векторное поле.

Теперь рассмотрим 1-форму на многообразии $M$. **1-форма** – это линейная вещественнозначная функция на векторах. Это означает, что 1-форма $\tilde{\sigma}$ в точке $P$ некоторое вещественное число, обозначаемое через $\tilde{\sigma}(\bar{v})$. (Волна над буквой обозначает 1-форму). В каждой точке $P$ многообразия $M$ можно построить линейное пространство, отличное от пространства 1-форм. Совокупность всех 1-форм в точке $P$ образует линейное векторное пространство – **кокасательное (дуальное, двойственное, сопряженное)** в точке $P \in M$, которое обозначается $T_p^*$. Аналогично векторам 1-формы образуют базис. Рассмотрим некоторую 1-форму $\tilde{\sigma}$, действующую на вектор $\bar{v}$:

$$\tilde{\sigma}(\bar{v}) = \tilde{\sigma}(v^b \bar{e}_b) = v^b \tilde{\sigma}(\bar{e}_b) = \tilde{\theta}^b(\bar{v})\tilde{\sigma}(\bar{e}_b)$$

Значения 1-форм определяются на векторах

$$\tilde{\sigma} = \sigma_\beta \tilde{\theta}^\beta$$

Значит, $\{\tilde{\theta}^b\}$ действительно является базисом. И этот базис будет двойственен базису пространства $T_p$.

### §1.3 Тензоры и внешние формы

Рассмотрим точку $P$ на $M$. Тензор типа $\binom{N}{N'}$ в точке определяется как линейная функция, аргументами которой служат $N$ 1-форм и $N'$ векторов, а значениями – вещественные числа. Это обобщенное определение 1-форм. **Компонентами тензора** являются его значения на базисных векторах и 1-формах. Как и в случае векторов и 1-форм, **тензорное поле** типа $\binom{N}{N'}$ – это правило, сопоставляющее каждой точке тензор типа $\binom{N}{N'}$ в этой точке. Свойство линейности тензоров распространяется и на тензорные поля.

Рассмотрим тензор $T$ типа $\binom{1}{1}$. Этот тензор предполагает наличие двух аргументов. Значит $T(\tilde{\sigma}, \overline{v})$ есть вещественное число, при фиксированном $\tilde{\sigma}$ мы получим 1-форму $T(\tilde{\sigma}, \_)$, поскольку сюда нужно подставить вектор, для того чтобы получить вещественное число; $T(\_, \overline{v})$ есть вектор. Таким образом, тензор типа $\binom{1}{1}$ можно рассматривать как линейную функцию на векторах со значениями в векторах, а также как линейную функцию на 1-формах, принимающую значение также в 1-формах. Это правило распространяется на любые тензоры.

Каждому тензору можно сопоставить матрицу его компонент.

1. Сложение и вычитание компонент тензора осуществляется по правилу сложения и вычитания матриц.
2. Умножение тензора на скаляр – умножение матрицы на скаляр.
3. Рассмотрим величины $\{v^i \sigma_j\}$, которые являются компонентами тензора типа $\binom{1}{1}$. Суммируя компоненты с учетом того, что $i = j$ получим, что $v^j \sigma_j$ – число, не зависящее от базиса, а именно значение $\tilde{\sigma}$ на $\overline{v}$, которое можно считать тензором типа $\binom{0}{0}$. Эта операция называется *сверткой*. Значит, *сверка* – это отображение прямого произведения касательного и кокасательного пространства в $\mathbb{R}$:
$$T_P \times T_P^* \longrightarrow \mathbb{R}.$$
4. Симметризация тензора.

Симметризация в компонентном виде означает следующее:

$$T_{(\mu,\nu,\ldots,\lambda)} = \frac{1}{n!}\big(T_{\mu,\nu,\ldots,\lambda} + T_{\nu,\mu,\ldots,\lambda} + \cdots\big)$$

$n$ – это количество индексов, участвующих в симметризации. Эта формула записана для нижних индексов, аналогично она выглядит и для верхних. Тензор называется *симметричным*, если при перестановке любых двух индексов он не меняет знак.

5. Альтернирование

Операция альтернирования выполняется по следующей формуле:

$$T_{[\mu,\nu,\ldots,\lambda]} = \frac{1}{n!}\bigl(T_{\mu,\nu,\ldots,\lambda} - T_{\nu,\mu,\ldots,\lambda} + \cdots\bigr)$$

$n$ – это количество индексов, участвующих в альтернировании. При альтернировании получаем тензор, который при любой четной перестановке нескольких индексов он не меняет знак, а при любой нечетной перестановке меняет его на противоположный. Такой тензор называется **кососимметричным.**

Тензор называется **антисимметричным**, если при перестановке любых двух индексов он меняет знак. **P-формой** ($p \geq 2$) (или **дифференциальной формой степени p**) называется антисимметричный тензор типа $\binom{0}{p}$. Множество всех p-форм образует векторное пространство относительно операции сложения и умножения на число:

$(\widetilde{\omega}_1 + \widetilde{\omega}_2)(\bar{v}) = \widetilde{\omega}_1(\bar{v}) + \widetilde{\omega}_2(\bar{v})$, размерность этого пространства $C_n^p$.

Теперь, пусть тензор $\check{T}$ - тензор типа $\binom{M}{N}$, а $\check{F}$ - тензор типа $\binom{M'}{N'}$. Тензор $\check{T} \otimes \check{F}$ - тензор типа $\binom{M+M'}{N+N'}$. Операция $\otimes$ называется **тензорным произведением**.

$$\check{T} \otimes \check{F}(\widetilde{\omega}_1, \ldots, \widetilde{\omega}_n, \widetilde{\sigma}_1, \ldots, \widetilde{\sigma}_n, \bar{v}_1, \ldots, \bar{v}_m, \bar{u}_1, \ldots, \bar{u}_m) =$$
$$= \check{F}(\widetilde{\omega}_1, \ldots, \widetilde{\omega}_n, \bar{v}_1, \ldots, \bar{v}_m) \cdot \check{T}(\widetilde{\sigma}_1, \ldots, \widetilde{\sigma}_n, \bar{u}_1, \ldots, \bar{u}_m)$$

**Внешним (косым) произведением** называется антисимметричное тензорное произведение дифференциальных форм. Пусть $\tilde{\alpha}$ и $\tilde{\beta}$ – 1-формы, тогда их внешним произведением будет являться 2-форма [2]:

$$\tilde{\alpha} \wedge \tilde{\beta} = \tilde{\alpha} \otimes \tilde{\beta} - \tilde{\beta} \otimes \tilde{\alpha} = 2!\,Alt\bigl(\tilde{\alpha} \otimes \tilde{\beta}\bigr)$$

Рассмотрим произвольный базис векторного пространства $\{\bar{e}_\alpha\}$ $\alpha = 1,..,n$, $\{\tilde{\theta}^\alpha\}$ - двойственный ему базис 1-форм и найдем значение 2-фомы $\tilde{\alpha}$ на векторах $\bar{u}$ и $\bar{v}$:

$$\tilde{\alpha}(\bar{u},\bar{v}) = \tilde{\alpha}(u^\mu \bar{e}_\mu, v^\nu \bar{e}_\nu) = u^\mu v^\nu \tilde{\alpha}(\bar{e}_\mu, \bar{e}_\nu) = u^\mu v^\nu \alpha_{\mu\nu} = u^\mu v^\nu \alpha_{[\mu\nu]} =$$
$$= u^{[\mu} v^{\nu]} \alpha_{\mu\nu} = \frac{1}{2!}(u^\mu v^\nu - u^\nu v^\mu) \alpha_{\mu\nu} =$$
$$= \frac{1}{2!}\left(\tilde{\theta}^\mu(\bar{u})\tilde{\theta}^\nu(\bar{v}) - \tilde{\theta}^\nu(\bar{u})\tilde{\theta}^\mu(\bar{v})\right)\alpha_{\mu\nu} =$$
$$= \frac{1}{2!}\alpha_{\mu\nu}(\tilde{\theta}^\mu \otimes \tilde{\theta}^\nu - \tilde{\theta}^\nu \otimes \tilde{\theta}^\mu)(\bar{u},\bar{v}) = \frac{1}{2!}\alpha_{\mu\nu}(\tilde{\theta}^\mu \wedge \tilde{\theta}^\nu)(\bar{u},\bar{v}).$$

Таким образом, получили разложение 2-формы $\tilde{\alpha}$ по базису векторного пространства всех 2-форм $\tilde{\theta}^\mu \wedge \tilde{\theta}^\nu$:

$$\tilde{\alpha} = \frac{1}{2!}\alpha_{\mu\nu}(\tilde{\theta}^\mu \wedge \tilde{\theta}^\nu), \quad \text{где } \alpha_{\mu\nu} = \tilde{\alpha}(\bar{e}_\alpha, \bar{e}_\beta) \quad (1.3.1)$$

Причем, число независимых базисных 2-форм (размерность векторного пространства 2-форм) равно числу ненулевых компонент $\alpha_{\mu\nu}$, т.е. $C_n^2$.

Следовательно, произвольная p-форма раскладывается по базису следующим образом: 
$$\tilde{\alpha} = \frac{1}{p!}\tilde{\alpha}_{\mu\nu\ldots\gamma}\tilde{\theta}^\mu \wedge \tilde{\theta}^\nu \wedge \ldots \wedge \tilde{\theta}^\gamma. \quad (1.3.2)$$

Найдем компактное представление (p+q)-формы, полученной применением операции внешнего произведения к произвольной p-форме и произвольной q-форме с учетом формулы (1.3.2).

Пусть $\quad \tilde{p} = \frac{1}{p!}p_{\alpha\ldots\beta}\tilde{\theta}^\alpha \wedge \ldots \wedge \tilde{\theta}^\beta$

$\tilde{q} = \frac{1}{q!}q_{\mu\ldots\nu}\tilde{\theta}^\mu \wedge \ldots \wedge \tilde{\theta}^\nu$

Тогда

$$\tilde{p} \wedge \tilde{q} = \frac{1}{(p+q)!}(\tilde{p} \wedge \tilde{q})_{\alpha\ldots\beta\mu\ldots\nu}\tilde{\theta}^\alpha \wedge \ldots \wedge \tilde{\theta}^\beta \wedge \tilde{\theta}^\mu \wedge \ldots \wedge \tilde{\theta}^\nu$$

С другой стороны

$$\tilde{p} \wedge \tilde{q} = \left(\frac{1}{p!}p_{\alpha\ldots\beta}\tilde{\theta}^\mu \wedge \ldots \wedge \tilde{\theta}^\nu\right) \wedge \left(\frac{1}{q!}q_{\mu\ldots\nu}\tilde{\theta}^\mu \wedge \ldots \wedge \tilde{\theta}^\nu\right) =$$
$$= \frac{1}{p!\,q!}\tilde{\theta}^\alpha \wedge \ldots \wedge \tilde{\theta}^\beta \wedge \tilde{\theta}^\mu \ldots \wedge \tilde{\theta}^\nu$$

Приравнивая два последних равенства получим:

$$(\tilde{p}\wedge\tilde{q})_{\alpha\ldots\beta\mu\ldots\nu} = \frac{(p+q)!}{p!\,q!} p_{[\alpha\ldots\beta} q_{\mu\ldots\nu]} = C_{p+q}^{p} p_{[\alpha\ldots\beta} q_{\mu\ldots\nu]}$$

Следовательно

$$\tilde{p}\wedge\tilde{q} = C_{p+q}^{p} p_{[\alpha\ldots\beta} q_{\mu\ldots\nu]} \tilde{\theta}^{\alpha}\wedge\ldots\wedge\tilde{\theta}^{\beta}\wedge\tilde{\theta}^{\mu}\wedge\ldots\wedge\tilde{\theta}^{\nu} \quad (1.3.3)$$

Операция внешнего произведения обладает следующими свойствами:

1) *косокоммутативность*: $\tilde{\alpha}\wedge\tilde{\beta} = (-1)^{pq}\tilde{\beta}\wedge\tilde{\alpha}$,
2) *дистрибутивность*: $(\lambda_1\tilde{\alpha}^1 + \lambda_2\tilde{\alpha}^2)\wedge\tilde{\beta} = \lambda_1\tilde{\alpha}^1\wedge\tilde{\beta} + \lambda_2\tilde{\alpha}^2\tilde{\beta},\ \lambda_1,\lambda_2 \in \mathbb{R}$,
3) *ассоциативность*: $(\tilde{\alpha}\wedge\tilde{\beta})\wedge\tilde{\gamma} = \tilde{\alpha}\wedge(\tilde{\beta}\wedge\tilde{\gamma})$,
4) $\tilde{\alpha}\wedge\tilde{\alpha} = \tilde{\alpha}\otimes\tilde{\alpha} - \tilde{\alpha}\otimes\tilde{\alpha} = 0$,

где $\tilde{\alpha}$ - p-форма, $\tilde{\beta}$ - q-форма, $\tilde{\gamma}$ - произвольная r-форма.

### §1.4 Метрический тензор

Рассмотрим одну из операций векторной алгебры, которая называется скалярным произведением векторов. Для этого возьмем два произвольных вектора $\bar{u}, \bar{v} \in T_P(M)$ и сопоставим им число по следующему правилу:

$$\breve{g}\colon \bar{u}, \bar{v} \to \bar{u}\cdot\bar{v} \in \mathbb{R}$$
$$\breve{g}(\bar{u}, \bar{v}) = \bar{u}\cdot\bar{v} \quad (1.4.1)$$

Из свойств скалярного произведения следует симметричность и билинейность $\breve{g}$:

$$\breve{g}(\bar{u}, \bar{v}) = \bar{u}\cdot\bar{v} = \bar{v}\cdot\bar{u} = \breve{g}(\bar{v}, \bar{u}),$$
$$\breve{g}(\alpha\bar{u}, \beta\bar{v}) = \alpha\bar{u}\cdot\beta\bar{v} = \alpha\beta(\bar{u}\cdot\bar{v}) = \alpha\beta\breve{g}(\bar{u}, \bar{v}).$$

(1.4.1) ведет себя как абстрактный тензор типа $\binom{0}{2}$, симметричный по своим элементам, он называется ***метрическим тензором***. Также можно ска-зать, что ***метрический тензор*** – это билинейная, симметричная, вещественнозначная форма на векторах [10].

***Компонентами метрического тензора*** называются его значения на базисных векторах, они представляют собой симметрическую матрицу типа $n \times n$. Если же матрица единичная, то метрический тензор называется ***евклидовой метрикой***.

$$g_{\mu\nu} = \begin{pmatrix} g_{11} & g_{12} & \cdots & g_{1n} \\ g_{21} & g_{22} & \cdots & g_{2n} \\ \vdots & \vdots & & \vdots \\ g_{n1} & g_{n2} & \cdots & g_{nn} \end{pmatrix}$$

$$g_{\mu\nu} = \breve{g}(\bar{e}_\mu, \bar{e}_\nu) = \bar{e}_\mu \cdot \bar{e}_\nu . \tag{1.4.2}$$

Положим $det\|g_{\mu\nu}\| \neq 0$, тогда существует обратная матрица метрического тензора, причем:

$$g_{\mu\nu} g^{\lambda\nu} = \delta_\mu^\lambda \tag{1.4.3}$$

С помощью метрики можно выполнять операцию поднятия и опускания индексов. Это возможность получить, например, из произвольного тензора $\breve{A}$ типа $\begin{pmatrix} 2 \\ 0 \end{pmatrix}$ тензор типа $\begin{pmatrix} 1 \\ 1 \end{pmatrix}$: $A^\mu{}_\nu = g_{\nu\lambda} A^{\mu\lambda}$. Этот тензор в свою очередь можно превратить в тензор типа $\begin{pmatrix} 0 \\ 2 \end{pmatrix}$: $A_{\tau\nu} = g_{\tau\alpha} A^\alpha{}_\nu = g_{\tau\alpha} g_{\nu\lambda} A^{\alpha\lambda}$. Таким образом, метрический тензор дает возможность отображать элементы касательного и кокасательного пространства друг в друга.

Пусть $M$ - многообразие класса $C^\infty$. ***Псевдориманова структура*** на $M$ есть тензорное поле $\breve{g}$ типа $\begin{pmatrix} 0 \\ 2 \end{pmatrix}$, удовлетворяющее условиям:

1. $\breve{g}(\bar{u}, \bar{v}) = \breve{g}(\bar{v}, \bar{u})$ для всех $\bar{u}, \bar{v} \in T_p(M)$
2. В каждой точке $P \in M$ значения $g_P$ поля $\breve{g}$ есть невырожденная билинейная форма $T_P \times T_P$.

***Псевдоримановым многообразием*** называется ***связное*** многообразие класса $C^\infty$ с псевдоримановой структурой. Если (и только если) форма $g_P$ по-ложительно определена для каждой точки $P \in M$, то можно говорить о римановой структуре и римановом многообразии[10].

Так как метрика на $M$ - это тензорное поле типа $\binom{0}{2}$, то теперь можно определить понятие расстояния на $M$. Получим выражение для **квадрата бесконечно малого вектора перемещения (квадрата интервала)** через метрический тензор:

$$dS^2 = d\bar{S} \cdot d\bar{S} = \breve{g}(d\bar{S}, d\bar{S}) = \breve{g}(dx^\alpha \bar{e}_\alpha, dx^\beta \bar{e}_\beta) = dx^\alpha dx^\beta \breve{g}(\bar{e}_\alpha, \bar{e}_\beta) =$$
$$= dx^\alpha dx^\beta g_{\alpha\beta}$$

$$dS^2 = dx^\alpha dx^\beta g_{\alpha\beta}. \qquad (1.4.4)$$

Здесь $d$ - символ бесконечно малой величины.

### §1.5 Связность и ковариантная производная

Рассмотрим метрический тензор на многообразии $M$. Этот тензор индуцируется координатным отображением на пространство $\mathbb{R}^n$. Он отображает элементы касательного $T_p(M)$ и кокасательного $T_p^*(M)$ пространств друг в друга, то есть поднимает и опускает индексы.

Часто для решения задач требуется сравнивать элементы разных тензорных полей (векторных полей и полей 1-форм) в разных точках. Для этого их нужно уметь переносить в одну точку, то есть надо знать правило параллельного переноса на многообразиях. Итак, на любом многообразии результат параллельного переноса, то есть переноса с сохранением углов, зависит от пути. Таким образом, для сравнения тензорных полей в разных точках многообразия нужно знать тот путь, вдоль которого происходит параллельный перенос данного объекта в исходную точку.

Зная изменение параметра кривой и изменение самого объекта, можно найти производную в точке по направлению. Но вид кривой и вид производной должны содержать информацию о пути, вдоль которого происходит дифференцирование. Всю эту информацию несет в себе **ковариантная производная**. Выясним, как находится ковариантная производная.

Рассмотрим векторное поле $W$, определенное в области U многообразия $M$, здесь же определена кривая $l$ с параметром $\lambda$. Если же поле $W$ задано в каждой точке кривой $l$, то можем определить его **ковариантную производную** вдоль $l$ в точке $P(\lambda = \lambda_0)$:

$$\nabla_{\bar{u}}\bar{w} = \lim_{\varepsilon \to 0} \frac{\bar{w}^*_{\lambda_0+\varepsilon}(\lambda_0) - \bar{w}(\lambda_0)}{\varepsilon}, \qquad (1.5.1)$$

где $\lambda_0$ – значение параметра в точке $P$, а вектор $\bar{w}^*_{\lambda_0+\varepsilon}(\lambda_0)$ получается параллельным переносом вектора $\bar{w}(\lambda_0+\varepsilon)$ назад в точку $P(\lambda_0)$, $\bar{u}$ показывает, вдоль какой кривой происходит дифференцирование.

Далее рассмотрим ковариантную производную от функции $f$ вдоль вектора $\bar{u}$. Это есть производная по направлению данной функции вдоль кривой с параметром $\lambda$ и касательным вектором $\bar{u}$:

$$\nabla_{\bar{u}} f = \frac{df}{d\lambda} = \bar{u}(f) = u^\alpha \frac{\partial f}{\partial x^\alpha} = u^\alpha \partial_\alpha f \qquad (1.5.2)$$

Операция ковариантного дифференцирования обладает следующими свойствами [1,2]:

$1^0$ ковариантная производная от скаляров совпадает с частной производной, то есть

$$\nabla_{\bar{u}} f = \frac{df}{d\lambda}, \qquad (1.5.3)$$

$2^0$ $\nabla_{\bar{u}}$ - дифференциальный оператор:

$$\nabla_{\bar{u}}(f\bar{w}) = f\nabla_{\bar{u}}\bar{w} + \bar{w}\nabla_{\bar{u}}f = f\nabla_{\bar{u}}\bar{w} + \bar{w}\frac{df}{d\lambda} = f\nabla_{\bar{u}}\bar{w} + \bar{w}\bar{u}(f), \qquad (1.5.4)$$

$3^0$ правило Лейбница позволяет распространить определение ковариантной производной на тензоры произвольного вида:

$$\nabla_{\bar{u}}(A \otimes B) = (\nabla_{\bar{u}} A) \otimes B + A \otimes (\nabla_{\bar{u}} B), \qquad (1.5.5)$$

$$\nabla_{\bar{u}}\langle \widetilde{\omega}, \bar{w} \rangle = \langle \nabla_{\bar{u}}\widetilde{\omega}, \bar{w} \rangle + \langle \widetilde{\omega}, \nabla_{\bar{u}}\bar{w} \rangle, \qquad (1.5.6)$$

$4^0$ понятие параллельного переноса вдоль кривой не зависит от выбора параметризации, то есть:

$$\nabla_{g\bar{u}}\bar{w} = g\nabla_{\bar{u}}\bar{w}, \qquad (1.5.7)$$

$5^0$ ковариантная производная аддитивна по разным направлениям:

$$(\nabla_{\bar{u}}\bar{w})_p + (\nabla_{\bar{v}}\bar{w})_p = (\nabla_{\bar{u}+\bar{v}}\bar{w})_p, \qquad (1.5.8)$$

Вычислим ковариантную производную от вектора:

$$\nabla_{\bar{u}}\bar{w} = \nabla_{u^\alpha \bar{e}_\alpha}(w^\beta \bar{e}_\beta) = u^\alpha\left((\nabla_{\bar{e}_\alpha} w^\beta)\bar{e}_\beta + w^\beta \nabla_{\bar{e}_\alpha}\bar{e}_\beta\right) =$$
$$= u^\alpha(\partial_\alpha w^\beta \bar{e}_\beta + w^\beta \nabla_{\bar{e}_\alpha}\bar{e}_\beta)$$

$\nabla_{\bar{e}_\alpha}\bar{e}_\beta$ – это дифференцирование базисного вектора $\bar{e}_\beta$ вдоль координатной линии с касательным вектором $\bar{e}_\alpha$. Разложим этот вектор по базису:

$$\nabla_{\bar{e}_\alpha}\bar{e}_\beta = \Gamma^\gamma_{\beta\alpha}\bar{e}_\gamma \qquad (1.5.9)$$

Величины $\Gamma^\gamma_{\beta\alpha}$ называются **коэффициентами связности**. Они несут информацию о том, какой объект дифференцируется, вдоль какой кривой многообразия и какой результат получается.

Таким образом: $\nabla_{\bar{u}}\bar{w} = (\partial_\alpha w^\gamma \bar{e}_\gamma + w^\beta \Gamma^\gamma_{\beta\alpha}\bar{e}_\gamma) = u^\alpha(\partial_\alpha w^\gamma + w^\beta \Gamma^\gamma_{\beta\alpha})\bar{e}_\gamma$.

$\nabla_\alpha w^\gamma = \partial_\alpha w^\gamma + w^\beta \Gamma^\gamma_{\beta\alpha}$ - компонента ковариантной производной от вектора.

Аналогично рассмотрим ковариантную производную от 1-формы, учитывая, что $\nabla_{\bar{e}_\beta}\tilde{\theta}^\alpha = \tilde{\Gamma}^\alpha_{\gamma\beta}\tilde{\theta}^\gamma$:

$$\nabla_{\bar{u}}\tilde{\omega} = \nabla_{u^\beta \bar{e}_\beta}(\omega_\alpha \tilde{\theta}^\alpha) = u^\beta\left((\nabla_{\bar{e}_\beta}\omega_\alpha)\tilde{\theta}^\alpha + \omega_\alpha \nabla_{\bar{e}_\beta}\tilde{\theta}^\alpha\right) =$$
$$= u^\beta\left((\partial_\beta \omega_\alpha)\tilde{\theta}^\alpha + \omega_\alpha \nabla_{\bar{e}_\beta}\tilde{\theta}^\alpha\right) = u^\beta\left((\partial_\beta \omega_\gamma)\tilde{\theta}^\gamma + \omega_\alpha \tilde{\Gamma}^\alpha_{\gamma\beta}\tilde{\theta}^\gamma\right) =$$
$$= u^\beta(\partial_\beta \omega_\gamma + \omega_\alpha \tilde{\Gamma}^\alpha_{\gamma\beta})\tilde{\theta}^\gamma$$

$\nabla_\beta \omega_\gamma = \partial_\beta \omega_\gamma + \omega_\alpha \tilde{\Gamma}^\alpha_{\gamma\beta}$ – компонента ковариантной производной от 1-формы.

**Связность Г** определяется коэффициентами связности, которые в свою очередь задают правило параллельного переноса по заданному пути на многообразиях.

Возникает вопрос, а как же связаны $\tilde{\Gamma}^{\gamma}_{\alpha\beta}$ и $\Gamma^{\gamma}_{\alpha\beta}$? Вычислим ковариантную производную вдоль кривой с касательным вектором $\bar{u}$ от свертки 1-формы и вектора двумя способами.

<u>1 способ.</u>

$$\nabla_{\bar{u}}\langle\tilde{\omega},\bar{v}\rangle = \nabla_{\bar{u}}(\omega_{\alpha}v^{\alpha}) = \nabla_{u^{\beta}\bar{e}_{\beta}}(\omega_{\alpha}v^{\alpha}) = u^{\beta}\nabla_{\bar{e}_{\beta}}(\omega_{\alpha}v^{\alpha}) = u^{\beta}\left(\left(\nabla_{\bar{e}_{\beta}}\omega_{\alpha}\right) + \omega_{\beta}\nabla_{\bar{e}_{\beta}}v^{\alpha}\right) = u^{\beta}\left((\partial_{\beta}\omega_{\alpha})v^{\alpha} + \omega_{\alpha}\partial_{\beta}v^{\alpha}\right) = u^{\beta}\left(v^{\alpha}\partial_{\beta}\omega_{\alpha} + \omega_{\alpha}\partial_{\beta}v^{\alpha}\right)$$

(1.5.10)

<u>2 способ.</u>

$$\nabla_{\bar{u}}\langle\tilde{\omega},\bar{v}\rangle = \langle\nabla_{\bar{u}}\tilde{\omega},\bar{v}\rangle + \langle\tilde{\omega},\nabla_{\bar{u}}\bar{v}\rangle = \langle u^{\beta}(\nabla_{\beta}\omega_{\gamma})\tilde{\theta}^{\gamma},\bar{v}\rangle + \langle\tilde{\omega},u^{\beta}(\nabla_{\beta}v^{\gamma})\bar{e}_{\gamma}\rangle =$$
$$= u^{\beta}\nabla_{\beta}\omega_{\gamma}\tilde{\theta}^{\gamma}(\bar{v}) + u^{\beta}\nabla_{\beta}v^{\gamma}\tilde{\omega}(\bar{e}_{\gamma}) = u^{\beta}\nabla_{\beta}\omega_{\gamma}v^{\gamma} + u^{\beta}\omega_{\gamma}\nabla_{\beta}v^{\gamma} = v^{\gamma}u^{\beta}\left(\partial_{\beta}\omega_{\gamma} + \tilde{\Gamma}^{\alpha}_{\gamma\beta}\omega_{\alpha}\right) + u^{\beta}\omega_{\gamma}\left(\partial_{\beta}v^{\gamma} + \Gamma^{\gamma}_{\alpha\beta}v^{\alpha}\right) = v^{\gamma}u^{\beta}\partial_{\beta}\omega_{\gamma} + v^{\gamma}u^{\beta}\tilde{\Gamma}^{\alpha}_{\gamma\beta}\omega_{\alpha} + u^{\beta}\omega_{\gamma}\partial_{\beta}v^{\gamma} + u^{\beta}\omega_{\gamma}\Gamma^{\gamma}_{\alpha\beta}v^{\alpha}$$

(1.5.11)

Приравняем (1.5.10) и (1.5.11):

$$v^{\alpha}u^{\beta}\partial_{\beta}\omega_{\alpha} + v^{\gamma}u^{\beta}\tilde{\Gamma}^{\alpha}_{\gamma\beta}\omega_{\alpha} + u^{\beta}\omega_{\alpha}\partial_{\beta}v^{\alpha} + u^{\beta}\omega_{\gamma}\Gamma^{\gamma}_{\alpha\beta}v^{\alpha} =$$
$$= u^{\beta}v^{\alpha}\partial_{\beta}\omega_{\alpha} + u^{\beta}\omega_{\alpha}\partial_{\beta}v^{\alpha}$$

$$v^{\gamma}u^{\beta}\tilde{\Gamma}^{\alpha}_{\gamma\beta}\omega_{\alpha} + u^{\beta}\omega_{\gamma}\Gamma^{\gamma}_{\alpha\beta}v^{\alpha} = 0$$

$$v^{\alpha}u^{\beta}\tilde{\Gamma}^{\gamma}_{\alpha\beta}\omega_{\gamma} + u^{\beta}\omega_{\gamma}\Gamma^{\gamma}_{\alpha\beta}v^{\alpha} = 0$$

$$u^{\beta}v^{\alpha}\omega_{\gamma}\left(\tilde{\Gamma}^{\gamma}_{\alpha\beta} + \Gamma^{\gamma}_{\alpha\beta}\right) = 0.$$

В силу произвольности u, v и $\omega$ имеем:

$$\tilde{\Gamma}^{\gamma}_{\alpha\beta} + \Gamma^{\gamma}_{\alpha\beta} = 0.$$

Таким образом, получим:  $\tilde{\Gamma}^{\gamma}_{\alpha\beta} = -\Gamma^{\gamma}_{\alpha\beta}$

Заметим, что данные коэффициенты связности не являются тензорами, так как закон, по которому происходит преобразование коэффициентов ($\Gamma^\lambda_{\mu\nu}$ в $\Gamma^{\lambda'}_{\nu'\mu'}$):

$$\Gamma^{\lambda'}_{\mu'\nu'} = \frac{\partial^2 x^\lambda}{\partial x^{\mu'} \partial x^{\nu'}} \frac{\partial x^{\lambda'}}{\partial x^\lambda} + \frac{\partial x^\mu}{\partial x^{\mu'}} \frac{\partial x^\nu}{\partial x^{\nu'}} \frac{\partial x^{\lambda'}}{\partial x^\lambda} \Gamma^\lambda_{\mu\nu}, \qquad (1.5.12)$$

не будет тензорным. Он совпал бы с тензорным, если бы в правой части не было первого слагаемого [7].

**Определение.** Метрический тензор и связность многообразия согласованы, если скалярное произведение двух любых векторов при их параллельном переносе вдоль произвольной кривой на многообразии сохраняется.

Условие согласованности метрики и связности на многообразии $M$ определяется по следующей формуле:

$$\nabla_\lambda g_{\mu\nu} = 0,$$

то есть, при параллельном переносе по замкнутому контуру длина вектора не меняется. Следствием из этого условия является формула:

$$\Gamma^\lambda_{\mu\nu} = \frac{1}{2} g^{\lambda\sigma} \Delta^{\alpha\beta\gamma}_{\sigma\nu\mu} (g_{\alpha\beta,\gamma} + T_{\alpha\beta\gamma}), \qquad (1.5.13)$$

где $\Delta^{\alpha\beta\gamma}_{\sigma\nu\mu} = \delta^\alpha_\sigma \delta^\beta_\nu \delta^\gamma_\mu - \delta^\alpha_\nu \delta^\beta_\mu \delta^\gamma_\sigma + \delta^\alpha_\mu \delta^\beta_\sigma \delta^\gamma_\nu$ - скобка Схоутэна.

Еще одно обозначение скобки Схоутэна – { }.

Вывод формулы для коэффициентов связности в пространстве с кручением приведен в **Приложении 1**.

### §1.6 Формы со значениями в векторных пространствах

***Формой со значениями в векторном пространстве*** или ***тензорнозначной*** формой называются формы, компонентами которых являются вектор или тензор произвольного ранга.

Простейшим примером тензорнозначной 1-формы является градиент вектора (ковариантная производная данного вектора в неконкретизированном направлении). Рассмотрим вектор $\bar{w}$ и выпишем градиент этого вектора

$$(\nabla_{\_}\bar{w})(\_). \tag{1.6.1}$$

В первый канал можно поместить вектор. Заметим, что (1.6.1.) обладает свойством линейности по первому каналу, то есть если $f, g$ - функции, $\bar{u}, \bar{v}$ – векторы, тога линейную комбинацию $f\bar{u} + g\bar{v}$ можно использовать в качестве аргумента в первом канале (1.6.1):

$$\nabla_{f\bar{u}+g\bar{v}}\bar{w}(\_) = f\nabla_{\bar{u}}\bar{w} + g\nabla_{\bar{v}}\bar{w} \tag{1.6.2}$$

Следовательно, объект (1.6.1) по первому каналу обладает свойствами формы. Если же туда поместить вектор, то мы получим иной объект, в который можно ввести 1-форму, то есть это будет тензор ранга $\binom{1}{0}$, иными словами вектор: $\nabla_{\bar{u}}\bar{w}(\_)$

Итак, объект (1.6.1) является тензорнозначной формой со значениями в векторном пространстве. Далее определим объект $\breve{T}(\_\_\_)$, который можно представить как тензорнозначную 0-форму. Пусть $\bar{w}(\_)$ – объект, компонентами которого являются вектора, то есть тензор ранга $\binom{1}{0}$, а $\breve{T}(\_\_\_)$ - объект, компонентами которого являются тензоры вида $\binom{3}{0}$. С точки зрения форм, это 0-формы. Поэтому они называются 0-формами, со значениями в векторном пространстве.

Вернемся к рассмотрению объекта (1.6.1). До того, как мы подействовали на него оператором $\nabla$ объект $\bar{w}(\_)$ являлся тензорнозначной 0-формой (то есть вектором). После применения $\nabla$ появился один канал формы. Следовательно (1.6.1) теперь 1-форма. Аналогично рассуждая, получим, что $\nabla_{\_}\breve{T}$ – тензорнозначная 1-форма.

### §1.7 Внешнее дифференцирование.
### Обобщенный внешний дифференциал

В данном параграфе определим дифференциальный оператор $\tilde{d}$, который сохраняет свойства форм и обратен к операции интегрирования в смысле формулы Ньютона-Лейбница.

Операция внешнего дифференцирования была введена Пуанкаре, и $\tilde{d}$ – дифференциальный оператор, который $p$-форму отображает в $(p+1)$-форму.

Пусть $M$ одномерное многообразие, тогда оператор $\tilde{d}$ превращает 0-форму $\tilde{\alpha}$ в 1-форму $\tilde{d}\tilde{\alpha}$.

Если расширить понятие оператора $\tilde{d}$ на формы старших степеней, то оператор дифференцирования будет обладать следующими свойствами:

Пусть $\tilde{\alpha}$ - p-форма, $\tilde{\beta}$, $\tilde{\gamma}$ - q-формы, тогда

$1^0$ $\tilde{d}\tilde{\alpha}$ – (p+1)-форма

$2^0$ $\tilde{d}(\tilde{\beta} + \tilde{\gamma}) = \tilde{d}\tilde{\beta} + \tilde{d}\tilde{\gamma}$ (линейность)

$3^0$ $\tilde{d}(\tilde{\alpha} \wedge \tilde{\beta}) = \tilde{d}\tilde{\alpha} \wedge \tilde{\beta} + (-1)^p \tilde{\alpha} \wedge \tilde{d}\tilde{\beta}$ (антидифференцирование)

$4^0$ $\tilde{d}(\tilde{d}\tilde{\alpha}) = \tilde{d}^2\tilde{\alpha} = 0$ (неимпотентность)

В (1.6.1) $\nabla$ повышает ранг векторнозначной 0-формы – вектора $\overline{w}$ на 1, т.е. (1.6.1) есть обобщение понятия внешней производной $\tilde{d}$ от векторнозначной 0-формы, т.е.

$$\nabla_{\_}\overline{w} = \mathcal{D}\overline{w} \qquad (1.7.1)$$

($\nabla_{\_}\overline{w}$ есть действие некоторого дифференциального оператора $\mathcal{D}$, увеличивающего ранг векторнозначной 0-формы на единицу).

Аналогично, обобщенная внешняя производная от тензорнозначной 0-формы $\tilde{T}(_\sim,_\sim,_\sim)$

$$\mathcal{D}\check{T} = \nabla_{-}\check{T} \qquad (1.7.2)$$

Обобщенная внешняя производная от скалярнозначной $p$-формы просто совпадает с обычной внешней производной $\tilde{d}$, которая обладает свойством антидифференцирования ($\tilde{\alpha}$ – $p$-форма, $\tilde{\beta}$ – произвольная $g$-форма):

$$\tilde{d}(\tilde{\alpha} \wedge \tilde{\beta}) = \tilde{d}\tilde{\alpha} \wedge \tilde{\beta} + (-1)^p \tilde{\alpha} \wedge \tilde{d}\tilde{\beta} \qquad (1.7.3)$$

Итак, мы определим понятие обобщенной внешней производной (обобщенного внешнего дифференциала) $\mathcal{D}$, как дифференциального оператора, обладающего следующими свойствами:

1) обобщенный внешний дифференциал от функции совпадает с ее градиентом:
$$\mathcal{D}f = \tilde{d}f$$
$\langle \mathcal{D}f, \bar{u} \rangle = \bar{u}(f)$, где $\bar{u}$ – произвольный вектор;

(производная по направлению вектора $\bar{u}$)

2) обобщенный внешний дифференциал от вектора (или векторнозначной 0-формы) совпадает с ковариантной производной данного вектора:
$$\mathcal{D}\bar{w} = \nabla_{-}\bar{w},$$
($\nabla_{-}\bar{w}$ – градиент вектора $\bar{w}$)

$\langle \mathcal{D}\bar{w}, \bar{u} \rangle = \nabla_{\bar{u}}\bar{w}$, где $\bar{u}$- произвольный вектор;

($\nabla_{\bar{u}}\bar{w}$ - ковариантная производная вдоль вектора $\bar{u}$)

3) обобщенный дифференциал от тензора $\check{S}$ (тензорнозначной 0-формы) совпадает с ковариантной производной данного тензора:
$$\mathcal{D}\check{S} = \nabla_{-}\check{S},$$
$\langle \mathcal{D}\check{S}, \bar{u} \rangle = \nabla_{\bar{u}}\check{S}$, где $\bar{u}$ – произвольный вектор;

4) обобщенный внешний дифференциал от скалярнозначной $p$-формы совпадает с обычным внешним дифференциалом от $p$-формы:
$$\mathcal{D}\tilde{\alpha} = \tilde{d}\tilde{\alpha},$$
где $\tilde{\alpha}$ – скалярнозначная $p$-форма

Условий (1) – (4) достаточно для того, чтобы распространить понятие обобщенного внешнего дифференциала на тензорнозначные $p$-формы $\check{\tilde{S}}$ со свойством, обобщающим (1.7.3):

$$\mathcal{D}\left(\check{S}\wedge\tilde{\beta}\right) = \mathcal{D}\check{S}\wedge\tilde{\beta} + (-1)^p \check{S}\wedge\mathcal{D}\tilde{\beta} = \mathcal{D}\check{S}\wedge\tilde{\beta} + (-1)^p \check{S}\wedge\tilde{d}\tilde{\beta} \quad (1.7.4)$$

($\tilde{\beta}$ – произвольная скалярнозначная $g$-форма).

## §1.8 Понятие кривизны многообразия

Кривизна характеризует изменение касательного вектора при переносе его по замкнутому контуру. При параллельном переносе вектора из одной точки в другую в пространстве с кривизной результаты получаются разные, если он совершается по разным путям. Это видно из **(рис. 5)**

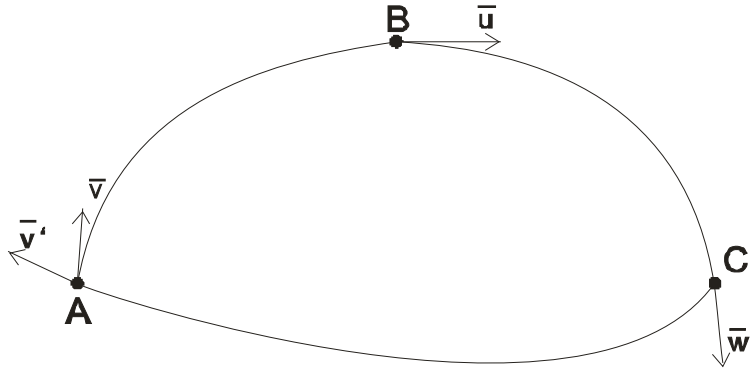

**Рис. 5**

Вектор $\bar{v}$ пройдет путь $AB$ в точке $B$ займет положение $\bar{u}$ (на всем пути $\bar{v}$ двигается, сохраняя угол). Аналогично двигаясь по пути $BC$ вектор $\bar{u}$ займет положение $\bar{w}$. И, наконец, пройдя путь $AC$, вектор $\bar{w}$ в точке $A$ будет находиться в положении $\bar{v}'$, что отлично от первоначального положения.

***Тензором кривизны Римана*** называется величина, выражаемая через компоненты связности [7]:

$$R^{\mu}{}_{\nu\sigma\lambda} = 2\partial_{[\alpha}\Gamma^{\mu}{}_{|\nu|\lambda]} + 2\Gamma^{\mu}{}_{\rho[\alpha}\Gamma^{\sigma}{}_{|\nu|\lambda]}. \quad (1.8.1)$$

Сверткой первого и третьего индексов тензора кривизны можно образовать симметрический тензор второго ранга:

$$R^{\mu}{}_{\nu\sigma\lambda}\delta^{\sigma}_{\mu} = R^{\mu}{}_{\nu\mu\lambda} = R_{\nu\lambda}, \quad (1.8.2)$$

который называется **тензором Риччи**. Если же свернуть тензор Риччи, то получим некоторый инвариант:

$$R = g^{\nu\lambda} R_{\nu\lambda} = R^{\nu}{}_{\nu}, \tag{1.8.3}$$

который называется **скалярной кривизной**.

Тензор кривизны обладает следующими свойствами в пространстве без кручения (связность симметрична):

$1^0$ $R^{\sigma}{}_{\mu\vartheta\lambda} = -R^{\sigma}{}_{\mu\lambda\vartheta}$,

$2^0$ $R_{\mu\vartheta\lambda\sigma} = R_{\lambda\sigma\mu\vartheta}$ (у тензора Римана можно менять местами пары индексов),

$3^0$ $R^{\sigma}{}_{[\mu\vartheta\lambda]} = 0$,

$4^0$ $R_{(\sigma\mu)\vartheta\lambda} = 0$ (тензор кривизны антисимметричен по первым 2-м индексам),

$5^0$ $\nabla_{[\rho} R^{\sigma}{}_{|\mu|\vartheta\lambda]} = 0$ (тождество Бианки).

### §1.9 Первое структурное уравнение Картана. 2-форма кручения

Выясним обобщенный внешний дифференциал (далее просто внешний дифференциал) от базисных векторов вдоль вектора $\bar{u}$:

$$\langle \mathcal{D}\bar{e}_\alpha, \bar{u} \rangle = \nabla_{\bar{u}} \bar{e}_\alpha = u^\gamma \nabla_{\bar{e}_\gamma} \bar{e}_\alpha = u^\gamma \Gamma^{\beta}{}_{\alpha\gamma} \bar{e}_\beta = \tilde{\theta}^\gamma(\bar{u}) \Gamma^{\beta}{}_{\alpha\gamma} \bar{e}_\beta = \bar{e}_\beta \Gamma^{\beta}{}_{\alpha\gamma} \tilde{\theta}^\gamma(\bar{u}) = \\ = \bar{e}_\beta \tilde{\Gamma}^{\beta}{}_{\alpha}(\bar{u}). \tag{1.9.1}$$

В силу произвольности вектора $\bar{u}$ его можно убрать:

$$\tilde{\Gamma}^{\beta}{}_{\alpha}(\_) = \tilde{\theta}^\gamma(\_) - \text{1-форма связности с компонентами } \Gamma^{\beta}{}_{\alpha\gamma}. \tag{1.9.2}$$

Кручение многообразия $M$ проще всего охарактеризовать в координатном базисе, а именно, его компонентами. В этом случае **компоненты тензора кручения** – это антисимметричная часть коэффициентов связности:

$$T^{\gamma}{}_{\alpha\beta} = -2\Gamma^{\gamma}{}_{[\alpha\beta]} = -\left(\Gamma^{\gamma}{}_{\alpha\beta} - \Gamma^{\gamma}{}_{\beta\alpha}\right). \tag{1.9.3}$$

Математический смысл понятия кручения, прояснит первое структурное уравнение Картана, вывод которого заключается в следующем.

Рассмотрим векторнозначную 1-форму, называемую **каноническая 1-формой кокасательного пространства**:

$$\bar{\bar{\delta}}(\sim, \_) = \bar{e}_\alpha(\sim) \otimes \tilde{\theta}^\alpha(\_). \tag{1.9.4}$$

(1.9.4) имеет два канала для ввода: один для ввода формы и один – для ввода вектора. (1.9.4) играет важную роль в теории гладких многообразий, т.к. он, действуя на объекты как линейный оператор, оставляет их без изменений.

Покажем это для произвольной 1-формы $\tilde{\omega} = \omega_\beta \tilde{\theta}^\beta$:

$$\bar{\bar{\delta}}(\tilde{\omega}, \_) = (\bar{e}_\alpha \otimes \tilde{\theta}^\alpha)(\tilde{\omega}) = \bar{e}_\alpha(\tilde{\omega})\tilde{\theta}^\alpha = \omega_\alpha \tilde{\theta}^\alpha = \tilde{\omega}.$$

Аналогично, для произвольного вектора $\bar{u} = u^\beta \bar{e}_\beta$:

$$\bar{\bar{\delta}}(\sim, \bar{u}) = (\bar{e}_\alpha \otimes \tilde{\theta}^\alpha)(\bar{u}) = \bar{e}_\alpha \tilde{\theta}^\alpha(\bar{u}) = \bar{e}_\alpha u^\alpha = \bar{u}.$$

$$\mathcal{D}\bar{\bar{\delta}} = \mathcal{D}(\bar{e}_\alpha \otimes \tilde{\theta}^\alpha) = (\mathcal{D}\bar{e}_\alpha)\wedge\tilde{\theta}^\alpha + (-1)^0 \bar{e}_\alpha \otimes \mathcal{D}\tilde{\theta}^\alpha = \nabla\bar{e}_\alpha\wedge\tilde{\theta}^\alpha + \bar{e}_\alpha \otimes \tilde{d}\tilde{\theta}^\alpha$$

$$\mathcal{D}\bar{e}_\alpha = \nabla\bar{e}_\alpha = \bar{e}_\beta \otimes \tilde{\Gamma}^\beta{}_\alpha \quad \text{(по (1.9.1))}$$

С учетом этого

$$\mathcal{D}\bar{\bar{\delta}} = \bar{e}_\beta \otimes (\tilde{\Gamma}^\beta{}_\alpha \wedge \tilde{\theta}^\alpha) + \bar{e}_\alpha \otimes \tilde{d}\tilde{\theta}^\alpha = \bar{e}_\alpha \otimes (\tilde{\Gamma}^\alpha{}_\beta \wedge \tilde{\theta}^\beta) + \bar{e}_\alpha \otimes \tilde{d}\tilde{\theta}^\alpha$$
$$= \bar{e}_\alpha \otimes (\tilde{d}\tilde{\theta}^\alpha + \tilde{\Gamma}^\alpha{}_\beta \wedge \tilde{\theta}^\beta)$$

Введем обозначение для ковариантной производной от 1-формы:

$$\mathcal{D}\tilde{\theta}^\alpha = \tilde{d}\tilde{\theta}^\alpha + \tilde{\Gamma}^\alpha{}_\beta \wedge \tilde{\theta}^\beta$$

Тогда $$\mathcal{D}\bar{\bar{\delta}} = \bar{e}_\alpha \otimes \mathcal{D}\tilde{\theta}^\alpha \tag{1.9.5}$$

Значение $\mathcal{D}$ от тензорнозначной 1-формы $\check{\tilde{S}}$ на векторах $\bar{u}$ и $\bar{v}$:

$$(\mathcal{D}\check{\tilde{S}})(\bar{u}, \bar{v}) = \nabla_{\bar{u}}\langle\check{\tilde{S}}, \bar{v}\rangle - \nabla_{\bar{v}}\langle\check{\tilde{S}}, \bar{u}\rangle - \langle\check{\tilde{S}}, [\bar{u}, \bar{v}]\rangle. \tag{1.9.6}$$

Вместо $\check{\tilde{S}}$ подставим векторнозначную 1-форму $\bar{\bar{\delta}}$:

$$\left(\mathcal{D}\bar{\bar{\delta}}\right)(\bar{u},\bar{v}) = \nabla_{\bar{u}}\langle\bar{\bar{\delta}},\bar{v}\rangle - \nabla_{\bar{v}}\langle\bar{\bar{\delta}},\bar{u}\rangle - \langle\bar{\bar{\delta}},[\bar{u},\bar{v}]\rangle = \nabla_{\bar{u}}\bar{v} - \nabla_{\bar{v}}\bar{u} - [\bar{u},\bar{v}]. \quad (1.9.7)$$

Вместо $\bar{u}$ и $\bar{v}$ подставим $\bar{e}_\alpha$ и, причем $\bar{e}_\alpha = \partial_\alpha$

$$\left(\mathcal{D}\bar{\bar{\delta}}\right)(\bar{e}_\alpha,\bar{e}_\beta) = \nabla_{\bar{e}_\alpha}\bar{e}_\beta - \nabla_{\bar{e}_\beta}\bar{e}_\alpha - [\bar{e}_\alpha,\bar{e}_\beta], \quad [\bar{e}_\alpha,\bar{e}_\beta] = [\partial_\alpha,\partial_\beta] = 0$$

$$\left(\mathcal{D}\bar{\bar{\delta}}\right)(\bar{e}_\alpha,\bar{e}_\beta) = \Gamma^\gamma{}_{\beta\alpha}\bar{e}_\gamma - \Gamma^\gamma{}_{\alpha\beta}\bar{e}_\gamma = \left(\Gamma^\gamma{}_{\beta\alpha} - \Gamma^\gamma{}_{\alpha\beta}\right)\bar{e}_\gamma = T^\gamma{}_{\alpha\beta}\bar{e}_\gamma = \breve{T}_{\alpha\beta}$$

$T^\gamma{}_{\alpha\beta}\bar{e}_\gamma = \breve{T}_{\alpha\beta}$ – компоненты тензорнозначной 2-формы кручения $\breve{\tilde{T}}$.

$$\breve{T}_{\alpha\beta} = -\breve{T}_{\beta\alpha}$$

$$\breve{\tilde{T}} = \frac{1}{2!}\breve{T}_{\alpha\beta}\tilde{\theta}^\alpha\wedge\tilde{\theta}^\beta = \frac{1}{2!}T^\gamma{}_{\alpha\beta}\bar{e}_\gamma\otimes(\tilde{d}x^\alpha\wedge\tilde{d}x^\beta)$$

Докажем равенство

$$\breve{\tilde{T}}(\bar{e}_\alpha,\bar{e}_\beta) = \frac{1}{2!}\breve{T}_{\mu\nu}(\tilde{\theta}^\mu\wedge\tilde{\theta}^\nu)(\bar{e}_\alpha,\bar{e}_\beta)$$

Учитывая то, что

$$(\tilde{\theta}^\mu\wedge\tilde{\theta}^\nu)(\bar{e}_\alpha,\bar{e}_\beta) = (\tilde{\theta}^\mu\otimes\tilde{\theta}^\nu - \tilde{\theta}^\nu\otimes\tilde{\theta}^\mu)(\bar{e}_\alpha,\bar{e}_\beta) = \delta^\mu_\alpha\delta^\nu_\beta - \delta^\nu_\alpha\delta^\mu_\beta = 2\delta^{[\mu}_\alpha\delta^{\nu]}_\beta$$

получим

$$\breve{\tilde{T}}(\bar{e}_\alpha,\bar{e}_\beta) = \frac{1}{2!}\breve{T}_{\mu\nu}2\delta^{[\mu}_\alpha\delta^{\nu]}_\beta = \breve{T}_{\alpha\beta}.$$

Таким образом

$$\left(\mathcal{D}\bar{\bar{\delta}}\right)(\bar{e}_\alpha,\bar{e}_\beta) = \breve{\tilde{T}}(\bar{e}_\alpha,\bar{e}_\beta),$$

$$\mathcal{D}\bar{\bar{\delta}} = \breve{\tilde{T}}. \quad (1.9.8)$$

Приравняем (1.9.5) и (1.9.8)

$$\breve{\tilde{T}} = \bar{e}_\alpha\otimes\mathcal{D}\tilde{\theta}^\alpha$$

$\breve{\tilde{T}} = \bar{e}_\gamma\otimes\left(\frac{1}{2!}T^\gamma{}_{\alpha\beta}\tilde{d}x^\alpha\wedge\tilde{d}x^\beta\right) = \bar{e}_\gamma\otimes\tilde{T}^\gamma$, где $\tilde{T}^\gamma$ – скалярнозначная 2-форма кручения

$$\breve{\tilde{T}} = \bar{e}_\alpha\otimes\tilde{T}^\alpha$$

$$\tilde{T}^\alpha = \mathcal{D}\tilde{\theta}^\alpha \text{ или } \tilde{T}^\alpha = \tilde{d}\tilde{\theta}^\alpha + \tilde{\Gamma}^\alpha{}_\beta \wedge \tilde{\theta}^\beta. \tag{1.9.9}$$

*(1.9.9) – первое структурное уравнение Картана.*

### §1.10 Второе структурное уравнение Картана. 2-форма кривизны

Для введения понятия кривизны связности рассмотрим двукратное применение обобщенного внешнего дифференциала к тензорнозначной 0-форме $\bar{w}$ (вектору)

$$\begin{aligned}
\mathcal{D}\bar{w} &= \mathcal{D}(\bar{e}_\alpha w^\alpha) = \mathcal{D}(\bar{e}_\alpha)w^\alpha + (-1)^0 \bar{e}_\alpha \otimes w^\alpha = (\nabla \bar{e}_\alpha)w^\alpha + \bar{e}_\alpha \otimes \tilde{d}w^\alpha \\
&= \bar{e}_\beta \otimes \tilde{\Gamma}^\beta{}_\alpha w^\alpha + \bar{e}_\alpha \otimes \tilde{d}w^\alpha = \bar{e}_\alpha \otimes \tilde{\Gamma}^\alpha{}_\beta w^\beta + \bar{e}_\alpha \otimes \tilde{d}w^\alpha \\
&= \bar{e}_\alpha \otimes (\tilde{d}w^\alpha + \tilde{\Gamma}^\alpha{}_\beta w^\beta)
\end{aligned}$$

$$\begin{aligned}
\mathcal{D}^2 \bar{w} &= \mathcal{D}(\mathcal{D}\bar{w}) = \mathcal{D}\left(\bar{e}_\alpha \otimes (\tilde{d}w^\alpha + \tilde{\Gamma}^\alpha{}_\beta w^\beta)\right) = \\
&= \mathcal{D}\bar{e}_\alpha \wedge (\tilde{d}w^\alpha + \tilde{\Gamma}^\alpha{}_\beta w^\beta) + (-1)^0 \bar{e}_\alpha \otimes \mathcal{D}(\tilde{d}w^\alpha + \tilde{\Gamma}^\alpha{}_\beta w^\beta) = \\
&= \nabla \bar{e}_\alpha \wedge (\tilde{d}w^\alpha + \tilde{\Gamma}^\alpha{}_\beta w^\beta) + \bar{e}_\alpha \otimes (\mathcal{D}\tilde{d}w^\alpha + (\mathcal{D}\tilde{\Gamma}^\alpha{}_\beta)w^\beta + (-1)^1 \tilde{\Gamma}^\alpha{}_\beta \wedge \mathcal{D}w^\beta) = \\
&= \bar{e}_\gamma \otimes (\tilde{\Gamma}^\gamma{}_\alpha \wedge \mathcal{D}w^\alpha) + \bar{e}_\gamma \otimes (\tilde{\Gamma}^\gamma{}_\alpha \wedge \tilde{\Gamma}^\alpha{}_\beta w^\beta) \\
&\quad + \bar{e}_\alpha \otimes \left(\underbrace{\tilde{d}^2 w^\alpha}_{0} + (\tilde{d}\tilde{\Gamma}^\alpha{}_\beta)w^\beta - \tilde{\Gamma}^\alpha{}_\beta \wedge \tilde{d}w^\beta\right) = \\
&= \bar{e}_\gamma \otimes (\tilde{\Gamma}^\gamma{}_\alpha \wedge \tilde{d}w^\alpha) + \bar{e}_\gamma \otimes \tilde{\Gamma}^\gamma{}_\alpha \wedge \tilde{\Gamma}^\alpha{}_\beta w^\beta + \bar{e}_\alpha \otimes (\tilde{d}\tilde{\Gamma}^\alpha{}_\beta)w^\beta \\
&\quad - \bar{e}_\alpha \otimes (\tilde{\Gamma}^\alpha{}_\beta \wedge \tilde{d}w^\beta) = \\
&= \underline{\bar{e}_\gamma \otimes (\tilde{\Gamma}^\gamma{}_\alpha \wedge \tilde{d}w^\alpha)} + \bar{e}_\alpha \otimes \tilde{\Gamma}^\alpha{}_\gamma \wedge \tilde{\Gamma}^\gamma{}_\beta w^\beta + \bar{e}_\alpha \otimes (\tilde{d}\tilde{\Gamma}^\alpha{}_\beta)w^\beta \\
&\quad - \underline{\bar{e}_\gamma \otimes (\tilde{\Gamma}^\gamma{}_\alpha \wedge \tilde{d}w^\alpha)} = \bar{e}_\alpha \otimes \left(\tilde{d}\tilde{\Gamma}^\alpha{}_\beta + \tilde{\Gamma}^\alpha{}_\gamma \wedge \tilde{\Gamma}^\gamma{}_\beta\right)w^\beta.
\end{aligned}$$

Обозначим скалярнозначную 2-форму в скобке через $\tilde{\mathcal{R}}^\alpha{}_\beta$ – 2-форма кривизны. Тогда

$$\mathcal{D}^2 \bar{w} = \bar{e}_\alpha \otimes \tilde{\mathcal{R}}^\alpha{}_\beta w^\beta, \tag{1.10.1}$$

где
$$\tilde{\mathcal{R}}^\alpha{}_\beta = \tilde{d}\tilde{\Gamma}^\alpha{}_\beta + \tilde{\Gamma}^\alpha{}_\gamma \wedge \tilde{\Gamma}^\gamma{}_\beta \tag{1.10.2}$$

*(1.10.2) – 2-ое структурное уравнение Картана.*

Получим в координатном базисе 2-ое структурное уравнение Картана.

$$\tilde{\mathcal{R}}^{\alpha}{}_{\beta} = \tilde{d}\tilde{\Gamma}^{\alpha}{}_{\beta} + \tilde{\Gamma}^{\alpha}{}_{\gamma}\wedge\tilde{\Gamma}^{\gamma}{}_{\beta}, \text{ базис } \{\tilde{d}x^{\alpha}\}$$

2-форма $\qquad \tilde{\mathcal{R}}^{\alpha}{}_{\beta} = \frac{1}{2!}\mathcal{R}^{\alpha}{}_{\beta\mu\nu}\tilde{d}x^{\mu}\wedge\tilde{d}x^{\nu},$

1-форма $\qquad \tilde{\Gamma}^{\alpha}{}_{\beta} = \Gamma^{\alpha}{}_{\beta\nu}\tilde{d}x^{\nu}(\_),$ тогда $\left(\tilde{d}^2 x^{\nu} = 0\right)$

$$\tilde{d}\left(\Gamma^{\alpha}{}_{\beta\nu}\tilde{d}x^{\nu}\right) = \tilde{d}\Gamma^{\alpha}{}_{\beta\nu}\wedge\tilde{d}x^{\nu} = \partial_{[\mu}\Gamma^{\alpha}{}_{|\beta|\nu]}\tilde{d}x^{\mu}\wedge\tilde{d}x^{\nu}$$

$$\tilde{\Gamma}^{\alpha}{}_{\gamma}\wedge\tilde{\Gamma}^{\gamma}{}_{\beta} = \Gamma^{\alpha}{}_{\gamma\mu}\tilde{d}x^{\mu}\wedge\Gamma^{\gamma}{}_{\beta\nu}\tilde{d}x^{\nu} = \Gamma^{\alpha}{}_{\gamma[\mu}\Gamma^{\gamma}{}_{|\beta|\nu]}\tilde{d}x^{\mu}\wedge\tilde{d}x^{\nu}$$

Левая часть: $\qquad \frac{1}{2!}\mathcal{R}^{\alpha}{}_{\beta\mu\nu}\tilde{d}x^{\mu}\wedge\tilde{d}x^{\nu}$

Правая часть: $\qquad \left(\partial_{[\mu}\Gamma^{\alpha}{}_{|\beta|\nu]} + \Gamma^{\alpha}{}_{\gamma[\mu}\Gamma^{\gamma}{}_{|\beta|\nu]}\right)\tilde{d}x^{\mu}\wedge\tilde{d}x^{\nu}$

$$\mathcal{R}^{\alpha}{}_{\beta\mu\nu} = 2\partial_{[\mu}\Gamma^{\alpha}{}_{|\beta|\nu]} + 2\Gamma^{\alpha}{}_{\gamma[\mu}\Gamma^{\gamma}{}_{|\beta|\nu]}. \tag{1.10.3}$$

Это известное выражение для тензора кривизны Римана.

### §1.11 Тождество Бианки для кручения

Найдем обобщенный внешний дифференциал $\mathcal{D}$ от тензорнозначной 2-формы кручения $\check{\tilde{T}} = \bar{e}_{\alpha}\otimes\tilde{T}^{\alpha}$.

$$\begin{aligned}
\mathcal{D}\check{\tilde{T}} &= \mathcal{D}\left(\bar{e}_{\alpha}\otimes\tilde{T}^{\alpha}\right) = \mathcal{D}\bar{e}_{\alpha}\wedge\tilde{T}^{\alpha} + (-1)^0 \bar{e}_{\alpha}\otimes\mathcal{D}\tilde{T}^{\alpha} = \nabla\bar{e}_{\alpha}\tilde{T}^{\alpha} + \bar{e}_{\alpha}\otimes\tilde{d}\tilde{T}^{\alpha} \\
&= \left(\bar{e}_{\beta}\otimes\tilde{\Gamma}^{\beta}{}_{\alpha}\right)\wedge\tilde{T}^{\alpha} + \bar{e}_{\alpha}\otimes\tilde{d}\left(\tilde{d}\tilde{\theta}^{\alpha} + \tilde{\Gamma}^{\alpha}{}_{\beta}\wedge\tilde{\theta}^{\beta}\right) = \\
&= \bar{e}_{\beta}\otimes\left(\tilde{\Gamma}^{\beta}{}_{\alpha}\wedge\tilde{T}^{\alpha}\right) + \bar{e}_{\alpha}\otimes\left(\underbrace{\tilde{d}^2\tilde{\theta}^{\alpha}}_{0} + \tilde{d}\tilde{\Gamma}^{\alpha}{}_{\beta}\wedge\tilde{\theta}^{\beta} + (-1)^1\tilde{\Gamma}^{\alpha}{}_{\beta}\wedge\tilde{d}\tilde{\theta}^{\beta}\right) = \\
&= \bar{e}_{\alpha}\otimes\left(\tilde{\Gamma}^{\alpha}{}_{\beta}\wedge\tilde{T}^{\beta}\right) + \bar{e}_{\alpha}\otimes\left(\tilde{d}\tilde{\Gamma}^{\alpha}{}_{\beta}\wedge\tilde{\theta}^{\beta} - \tilde{\Gamma}^{\alpha}{}_{\beta}\wedge\tilde{d}\tilde{\theta}^{\beta}\right) = \\
&= \bar{e}_{\alpha}\otimes\left(\tilde{\Gamma}^{\alpha}{}_{\beta}\wedge\tilde{T}^{\beta}\right) + \bar{e}_{\alpha}\otimes\left(\tilde{d}\tilde{\Gamma}^{\alpha}{}_{\beta}\wedge\tilde{\theta}^{\beta} - \tilde{\Gamma}^{\alpha}{}_{\beta}\wedge\left(\tilde{T}^{\beta} - \tilde{\Gamma}^{\beta}{}_{\gamma}\wedge\tilde{\theta}^{\gamma}\right)\right) = \\
&= \bar{e}_{\alpha}\otimes\left(\underline{\tilde{\Gamma}^{\alpha}{}_{\beta}\wedge\tilde{T}^{\beta}} + \tilde{d}\tilde{\Gamma}^{\alpha}{}_{\beta}\wedge\tilde{\theta}^{\beta} - \underline{\tilde{\Gamma}^{\alpha}{}_{\beta}\wedge\tilde{T}^{\beta}} + \tilde{\Gamma}^{\alpha}{}_{\beta}\wedge\tilde{\Gamma}^{\beta}{}_{\gamma}\wedge\tilde{\theta}^{\gamma}\right) = \\
&= \bar{e}_{\alpha}\otimes\left(\tilde{d}\tilde{\Gamma}^{\alpha}{}_{\beta}\wedge\tilde{\theta}^{\beta} + \wedge\tilde{\theta}^{\beta}\right) = \bar{e}_{\alpha}\otimes\left(\tilde{d}\tilde{\Gamma}^{\alpha}{}_{\beta} + \tilde{\Gamma}^{\alpha}{}_{\gamma}\wedge\tilde{\Gamma}^{\gamma}{}_{\beta}\right)\wedge\tilde{\theta}^{\beta} = \\
&= \bar{e}_{\alpha}\otimes\left(\tilde{\mathcal{R}}^{\alpha}{}_{\beta}\wedge\tilde{\theta}^{\beta}\right) = \bar{e}_{\alpha}\otimes\tilde{\theta}^{\beta}\wedge\tilde{\mathcal{R}}^{\alpha}{}_{\beta}
\end{aligned}$$

$$\mathcal{D}\check{\tilde{T}} = \bar{e}_{\alpha}\otimes\left(\tilde{\theta}^{\beta}\wedge\tilde{\mathcal{R}}^{\alpha}{}_{\beta}\right) \tag{1.11.1}$$

***(1.11.1) – 3-е структурное уравнение Картана или тождество Бианки для кручения.***

Получим формулу (1.11.1) в координатном базисе.

$$\mathcal{D}\tilde{\breve{T}} = \mathcal{D}(\bar{e}_\alpha \otimes \tilde{T}^\alpha) = \mathcal{D}\bar{e}_\alpha \wedge \tilde{T}^\alpha + \bar{e}_\alpha \otimes \tilde{d}\tilde{T}^\alpha =$$
$$= \nabla\bar{e}_\alpha \wedge \frac{1}{2}T^\alpha{}_{\beta\gamma}\tilde{d}x^\beta \wedge \tilde{d}x^\gamma + \bar{e}_\alpha \otimes \tilde{d}\left(\frac{1}{2}T^\alpha{}_{\beta\gamma}\tilde{d}x^\beta \wedge \tilde{d}x^\gamma\right) =$$
$$= \frac{1}{2}\left(\bar{e}_\lambda \otimes \underbrace{\tilde{\Gamma}^\lambda{}_\alpha}_{\Gamma^\lambda{}_{\alpha\nu}\tilde{d}x^\nu} \wedge \tilde{d}x^\beta \wedge \tilde{d}x^\gamma T^\alpha{}_{\beta\gamma} + \bar{e}_\alpha \otimes \partial_\lambda T^\alpha{}_{\beta\gamma}\tilde{d}x^\lambda \wedge \tilde{d}x^\beta \wedge \tilde{d}x^\gamma\right) =$$
$$= \frac{1}{2}\left(T^\alpha{}_{\beta\gamma}\Gamma^\lambda{}_{\alpha\nu}\bar{e}_\lambda \otimes \tilde{d}x^\nu \wedge \tilde{d}x^\beta \wedge \tilde{d}x^\gamma + \partial_\lambda T^\alpha{}_{\beta\gamma}\bar{e}_\alpha \otimes \tilde{d}x^\lambda \wedge \tilde{d}x^\beta \wedge \tilde{d}x^\gamma\right) =$$
$$= \frac{1}{2}\left(T^\nu{}_{\beta\gamma}\Gamma^\alpha{}_{\nu\lambda}\bar{e}_\alpha \otimes \tilde{d}x^\lambda \wedge \tilde{d}x^\beta \wedge \tilde{d}x^\gamma + \partial_\lambda T^\alpha{}_{\beta\gamma}\bar{e}_\alpha \otimes \tilde{d}x^\lambda \wedge \tilde{d}x^\beta \wedge \tilde{d}x^\gamma\right) =$$
$$= \frac{1}{2}\bar{e}_\alpha \otimes (\partial_\lambda T^\alpha{}_{\beta\gamma} + \Gamma^\alpha{}_{\nu\lambda}T^\nu{}_{\beta\gamma} - \Gamma^\nu{}_{\beta\lambda} - \Gamma^\nu{}_{\gamma\lambda} + \Gamma^\nu{}_{[\beta\lambda]}T^\alpha{}_{\nu\gamma}$$
$$+ \Gamma^\nu{}_{[\gamma\lambda]}T^\alpha{}_{\beta\nu})\tilde{d}x^\lambda \wedge \tilde{d}x^\beta \wedge \tilde{d}x^\gamma =$$
$$= \frac{1}{2}\bar{e}_\alpha \otimes \left(\nabla_\lambda T^\alpha{}_{\beta\gamma} - \left(-\frac{1}{2}\right)T^\nu{}_{\lambda\beta}T^\alpha{}_{\gamma\nu} - \frac{1}{2}2T^\nu{}_{\lambda\beta}T^\alpha{}_{\gamma\nu}\right)\tilde{d}x^\lambda \wedge \tilde{d}x^\beta \wedge \tilde{d}x^\gamma$$
$$= \frac{1}{2}\bar{e}_\alpha \otimes (\nabla_\lambda T^\alpha{}_{\beta\gamma} - T^\nu{}_{\lambda\beta}T^\alpha{}_{\gamma\nu})\tilde{d}x^\lambda \wedge \tilde{d}x^\beta \wedge \tilde{d}x^\gamma$$

Компоненты левой части: $\frac{1}{2}\left(\nabla_{[\lambda}T^\alpha{}_{\beta\gamma]} - T^\nu{}_{[\lambda\beta}T^\alpha{}_{\gamma]\nu}\right)$.

$$\bar{e}_\alpha \otimes (\tilde{\theta}^\beta \wedge \tilde{\mathcal{R}}^\alpha{}_\beta) = \bar{e}_\alpha \otimes (\tilde{\theta}^\beta \wedge \mathcal{R}^\alpha{}_{\beta\mu\nu}\tilde{d}x^\mu \wedge \tilde{d}x^\nu) = \bar{e}_\alpha \otimes (\tilde{\theta}^\lambda \wedge \mathcal{R}^\alpha{}_{\lambda\beta\gamma}\tilde{d}x^\beta \wedge \tilde{d}x^\gamma) =$$
$$= \frac{1}{2}\bar{e}_\alpha \otimes (\mathcal{R}^\alpha{}_{[\lambda\beta\gamma]})\tilde{d}x^\lambda \wedge \tilde{d}x^\beta \wedge \tilde{d}x^\gamma$$

Компоненты правой части: $\frac{1}{2}\mathcal{R}^\alpha{}_{[\lambda\beta\gamma]}$.

Тождество Бианки для кручения в компонентном виде имеет вид:

$$\mathcal{R}^\alpha{}_{[\lambda\beta\gamma]} = \nabla_{[\lambda}T^\alpha{}_{\beta\gamma]} - T^\nu{}_{[\lambda\beta}T^\alpha{}_{\gamma]\nu} \qquad (1.11.2)$$

в $V_4$ $\quad T^\alpha{}_{\beta\gamma} = 0 \implies \quad \mathcal{R}^\alpha{}_{[\lambda\beta\gamma]} = 0.$ $\qquad (1.11.3)$

# Глава II

## Связность и кривизна многообразий с волновой метрикой

### §2.1 Уравнение Эйнштейна и волновая метрика

Неевклидовы геометрии играют большую роль в прикладных вопросах современной теории естествознания. На основе неевклидовой геометрии строится теория гравитации. Современная теория гравитации – общая теория относительности, которая придала геометрический смысл гравитации, связав ее с римановой кривизной [6].

В основе общей теории относительности Эйнштейн положил принцип эквивалентности, суть которого заключается в локальной эквивалентности гравитационного поля и сил инерции, возникающих в неинерциальных система отсчета, переход к которым может быть осуществлен введением криволинейной системы координат. Величины, которые определяют в каждой данной криволинейной системе координат свойства геометрии, устанавливают, как известно, метрику пространства-времени $g_{\mu\nu}$. В общей теории относительности $g_{\mu\nu}$ отождествляется с потенциалами гравитационного поля, а коэффициенты связности $\Gamma^{\lambda}_{\mu\nu}$, определяющие параллельный перенос вдоль кривых в данной геометрии, отождествляются с силой, причем эти коэффициенты симметричны по первым двум индексам, то есть $\Gamma^{\mu}_{\lambda\nu} = \Gamma^{\lambda}_{\mu\nu}$. При этом кривизна пространства-времени $\mathcal{R}^{\lambda}_{\nu\mu}$ отлична от нуля и выражается через $\Gamma^{\lambda}_{\mu\nu}$ по формуле (1.8.1). Основные динамические уравнения ОТО были получены независимо Гильбертом и Эйнштейном. При этом Гильберт получил уравнения на основе вариационного исчисления для системы «гравитационное поле - материя», а Эйнштейн при получении уравнений гравитационного поля руководствовался тем принципом, что тензор энергии-импульса материи $T_{\mu\nu}$ должен порождать отличную от нуля кривизну пространства-времени, а кривизна пространства-времени должна определять движение материи. Эти уравнения имеют следующий вид [8]:

$$G_{\mu\nu} = R_{\mu\nu} - \frac{1}{2}g_{\mu\nu}R = \varkappa T_{\mu\nu},$$

где $G_{\mu\nu}$- бином Эйнштейна, построенный на основе $R_{\mu\nu}$ – тензора Риччи и $R$ - скаляра, определенного в §8, а $\varkappa$ - коэффициент пропорциональности, который называется гравитационной постоянной Эйнштейна.

Рассмотрим вид уравнения Эйнштейна в пустоте. «Пустое» здесь означает отсутствие материи и каких-либо физических полей, за исключением самого гравитационного поля. Гравитационное поле не нарушает пустоты, все остальные поля нарушают. Итак, в пустоте уравнение Эйнштейна в силу равенства нулю тензора энергии-импульса материи принимает вид:

$$R_{\mu\nu} - \frac{1}{2}g_{\mu\nu}R = 0. \tag{2.1.1}$$

Свернем (2.1.1) по $\mu$ и $\nu$. В результате получим:

$$R = 0. \tag{2.1.2}$$

Таким образом в силу (2.1.2) уравнение Эйнштейна в пустоте имеет вид:

$$R_{\mu\nu} = 0. \tag{2.1.3}$$

Напомним следующее определение:

***Псевдоевклидовой метрикой*** [5] сигнатуры $(p, q)$ на гладком многообразии $M$ размерности $n = p + q$ называется такая гладкая симметричная дифференциальная форма $g$ на $M$, что для любой точки $x \in M$ форма $g_x$ на $T_x(M)$ невырождена и имеет сигнатуру $(p, g)$. Пара $(M, g)$ называется ***псевдоримановым многообразием***.

Если $g = 0$ (то есть $g_x$ положительно определена), форма $g$ называется ***римановой метрикой***, а пара $(M, g)$ - ***римановым многообразием***.

Если $p = 1$, а $g > 0$, то форма $g$ называется ***лоренцевой метрикой***, а пара $(M, g)$ - ***лоренцевым многообразием***.

В дальнейшем мы будем рассматривать псевдориманово многообразие с лоренцевой сигнатурой.

Асимптотическое поведение гравитационных полей, создаваемых изолированными источниками, как известно, имеет большое сходство с поведением плоских электромагнитных волн в пространстве-времени Минковского. Это стало мотивировкой для Бонди, Пирани и Робинсона [11] дать строго групповое определение понятия плоских гравитационных волн в пустом пространстве для метрического поля, удовлетворяющего двум постулатам:

- ✓ Поле одинаково в любой точке волнового фронта.
- ✓ Метрический тензор пространства-времени подобно векторному потенциалу плоских электромагнитных волн допускает определенную группу симметрии.

Накладывается условие, что эти волны обладают степенью симметрий, аналогичной той, которой обладает электромагнитная волна в плоском пространстве-времени. Симметрия плоских электромагнитных волновых полей настолько велика, что требование такой симметрии для гравитационных полей автоматически влечет их волновой характер. Волновой характер проявляется во внутренних изотропных свойствах этих метрик. Это не «координатный вектор», чья метрика может быть преобразована к статической форме с помощью координатного преобразования. Изотропность становится видимой из физических свойств метрик.

Если рассматривать плоскую электромагнитную волну, распространяющуюся в положительном направлении оси $z$, то очевидные свойства симметрий такой волны заключаются в том, что существует 3-параметрическая группа преобразования пространства-времени Минковского в себя, оставляющая электромагнитное поле неизменным. Это движения-трансляции в направлениях $y$ и $z$ вдоль изотропной 3-поверхности $t - x = const$. Кроме этих очевидных симметрий существуют дополнительные симметрии плоской электромагнитной волны, которые менее очевидны, так как они являются внутренними 4-мерными. Эти симметрии описываются 2-параметрической группой движений, переводящей изотропные 3-поверхности $t - x = const$ в себя. Эти движения называются «*изотропными вращениями*».

Плоские волны могут быть интерпретированы как гравитационные поля на больших расстояниях от излучающих тел.

Как известно [11], проблему гравитационных волн удобнее рассматривать в специальным образом выбранном координатном базисе [12], образованном двумя нулевыми (изотропными) векторами $\bar{e}_0 = \partial_v$, $\bar{e}_1 = \partial_u$ и двумя пространственно-подобными $\bar{e}_2 = \partial_x$, $\bar{e}_3 = \partial_y$, причем вектор $\bar{e}_0$ ковариантно постоянен и направлен вдоль волнового луча, а координаты $x$ и $y$ параметризуют волновую поверхность $(u, v) = const$. Метрический тензор в этом базисе 1-форм равен

$$g_{ab} = \begin{pmatrix} 0 & 1 & 0 & 0 \\ 1 & 0 & 0 & 0 \\ 0 & 0 & -1 & 0 \\ 0 & 0 & 0 & -1 \end{pmatrix} \qquad (2.1.4)$$

Плоская волна метрики определяется как частный случай метрики плоско-фронтовых гравитационных волн с параллельными лучами (pp-волны), имеющий в указанном базисе вид [11]:

$$dS^2 = 2H(x, y, u)du^2 + 2dudv - dx^2 - dy^2, \qquad (2.1.5)$$

где координата $u$ имеет смысл запаздывающего времени и интерпретируется как фаза волны. Поэтому (2.1.5) – это волновая метрика. Метрике (2.1.5) соответствует базис из 1-форм:

$$\tilde{\theta}^0 = H\tilde{d}u + \tilde{d}v, \quad \tilde{\theta}^1 = \tilde{d}u, \quad \tilde{\theta}^2 = \tilde{d}x, \quad \tilde{\theta}^3 = \tilde{d}y. \qquad (2.1.6)$$

в котором эта метрика имеет вид

$$g = g_{ab}\tilde{\theta}^a \otimes \tilde{\theta}^b$$

Найдем условие, при котором метрика (2.1.5) является решением уравнения Эйнштейна в пустоте. Для этого вычислим геометрические характеристики многообразия без кручения, наделенного метрикой (2.1.5).

## §2.2 1-формы связности и коэффициенты связности
## для волновой метрики.

Получим вид 1-формы связности для волновой метрики (формула (2.1.5)) в случае, когда кручение $M$ равно нулю $T^a{}_{bc} = 0$. Для этого используем формулу определяющую связность в следующем виде [13]:

$$\tilde{\Gamma}_{ab} = -\frac{1}{2}\bar{e}_{\{c}\rfloor(\bar{e}_a\rfloor\tilde{d}\tilde{\theta}_{b\}})\tilde{\theta}^c. \tag{2.2.1}$$

Формула (2.2.1) получена в [13], где скобки { } называются *скобками Схоутэна*, и в компонентном виде совпадает с формулой (1.5.13). Используя описанный в параграфе §1 данной главы базис 1-форм (2.1.6), найдём внешний дифференциал от 1-формы базиса (2.1.6), пользуясь свойствами этого дифференциала:

$$\tilde{d}\tilde{\theta}^0 = -H_x\tilde{\theta}^1\wedge\tilde{\theta}^2 - H_y\tilde{\theta}^1\wedge\tilde{\theta}^3 = -H_x\tilde{d}u\wedge\tilde{d}x - H_y\tilde{d}u\wedge\tilde{d}y\,,$$
$$\tilde{d}\tilde{\theta}^1 = 0\,, \quad \tilde{d}\tilde{\theta}^2 = 0\,, \quad \tilde{d}\tilde{\theta}^3 = 0, \tag{2.2.2}$$

где введены следующие обозначения:

$H_x$ - производная $\frac{\partial H}{\partial x}$, $H_y$ – производная $\frac{\partial H}{\partial y}$.

После раскрытия скобки Схоутэна формула (2.2.1) принимает вид:

$$\tilde{\Gamma}_{ab} = -\frac{1}{2}\left[g_{bd}(\bar{e}_a\rfloor\tilde{d}\tilde{\theta}^d) - g_{ad}(\bar{e}_b\rfloor\tilde{d}\tilde{\theta}^d) - g_{cd}\bar{e}_a\rfloor(\bar{e}_b\rfloor\tilde{d}\tilde{\theta}^d)\tilde{\theta}^c\right] \tag{2.2.3}$$

Используя формулу (2.2.3) получим:

$$\tilde{\Gamma}_{00} = -\frac{1}{2}\left[g_{0d}(\bar{e}_0\rfloor\tilde{d}\tilde{\theta}^d) - g_{0d}(\bar{e}_0\rfloor\tilde{d}\tilde{\theta}^d) - g_{cd}\bar{e}_0\rfloor(\bar{e}_0\rfloor\tilde{d}\tilde{\theta}^d)\tilde{\theta}^c\right] =$$
$$= -\frac{1}{2}\left[g_{0d}(\bar{e}_0\rfloor\tilde{d}\tilde{\theta}^d) - g_{0d}(\bar{e}_0\rfloor\tilde{d}\tilde{\theta}^d) - g_{0d}\bar{e}_0\rfloor(\bar{e}_0\rfloor\tilde{d}\tilde{\theta}^d)\tilde{\theta}^0\right.$$
$$\left. - g_{1d}\bar{e}_0\rfloor(\bar{e}_0\rfloor\tilde{d}\tilde{\theta}^d)\tilde{\theta}^1 - g_{2d}\bar{e}_0\rfloor(\bar{e}_0\rfloor\tilde{d}\tilde{\theta}^d)\tilde{\theta}^2 - g_{3d}\bar{e}_0\rfloor(\bar{e}_0\rfloor\tilde{d}\tilde{\theta}^d)\tilde{\theta}^3\right] =$$

$$
\begin{aligned}
= -\frac{1}{2}\Bigg[ & \bigg( \underset{0}{g_{00}} (\bar{e}_0 \rfloor \tilde{d}\tilde{\theta}^0) - \underset{0}{g_{00}} (\bar{e}_0 \rfloor \tilde{d}\tilde{\theta}^0) - \underset{0}{g_{00}} \bar{e}_0 \rfloor (\bar{e}_0 \rfloor \tilde{d}\tilde{\theta}^0)\tilde{\theta}^0 \\
& - \underset{1}{g_{10}} \bar{e}_0 \rfloor (\bar{e}_0 \rfloor \tilde{d}\tilde{\theta}^0)\tilde{\theta}^1 - \underset{0}{g_{20}} \bar{e}_0 \rfloor (\bar{e}_0 \rfloor \tilde{d}\tilde{\theta}^0)\tilde{\theta}^2 - \underset{0}{g_{30}} \bar{e}_0 \rfloor (\bar{e}_0 \rfloor \tilde{d}\tilde{\theta}^d)\tilde{\theta}^3 \bigg) \\
& + \bigg( \underset{1}{g_{01}} (\bar{e}_0 \rfloor \tilde{d}\tilde{\theta}^1) - \underset{1}{g_{01}} (\bar{e}_0 \rfloor \tilde{d}\tilde{\theta}^1) - \underset{1}{g_{01}} \bar{e}_0 \rfloor (\bar{e}_0 \rfloor \tilde{d}\tilde{\theta}^1)\tilde{\theta}^0 \\
& - \underset{0}{g_{11}} \bar{e}_0 \rfloor (\bar{e}_0 \rfloor \tilde{d}\tilde{\theta}^1)\tilde{\theta}^1 - \underset{0}{g_{21}} \bar{e}_0 \rfloor (\bar{e}_0 \rfloor \tilde{d}\tilde{\theta}^1)\tilde{\theta}^2 - \underset{0}{g_{31}} \bar{e}_0 \rfloor (\bar{e}_0 \rfloor \tilde{d}\tilde{\theta}^1)\tilde{\theta}^3 \bigg) \\
& + \bigg( \underset{0}{g_{02}} (\bar{e}_0 \rfloor \tilde{d}\tilde{\theta}^2) - \underset{0}{g_{02}} (\bar{e}_0 \rfloor \tilde{d}\tilde{\theta}^2) - \underset{0}{g_{02}} \bar{e}_0 \rfloor (\bar{e}_0 \rfloor \tilde{d}\tilde{\theta}^2)\tilde{\theta}^0 \\
& - \underset{0}{g_{12}} \bar{e}_0 \rfloor (\bar{e}_0 \rfloor \tilde{d}\tilde{\theta}^2)\tilde{\theta}^1 - \underset{-1}{g_{22}} \bar{e}_0 \rfloor (\bar{e}_0 \rfloor \tilde{d}\tilde{\theta}^2)\tilde{\theta}^2 - \underset{0}{g_{32}} \bar{e}_0 \rfloor (\bar{e}_0 \rfloor \tilde{d}\tilde{\theta}^2)\tilde{\theta}^3 \bigg) \\
& + \bigg( \underset{0}{g_{03}} (\bar{e}_0 \rfloor \tilde{d}\tilde{\theta}^3) - \underset{0}{g_{03}} (\bar{e}_0 \rfloor \tilde{d}\tilde{\theta}^3) - \underset{0}{g_{03}} \bar{e}_0 \rfloor (\bar{e}_0 \rfloor \tilde{d}\tilde{\theta}^3)\tilde{\theta}^0 \\
& - \underset{0}{g_{13}} \bar{e}_0 \rfloor (\bar{e}_0 \rfloor \tilde{d}\tilde{\theta}^3)\tilde{\theta}^1 - \underset{0}{g_{23}} \bar{e}_0 \rfloor (\bar{e}_0 \rfloor \tilde{d}\tilde{\theta}^3)\tilde{\theta}^2 - \underset{-1}{g_{33}} \bar{e}_0 \rfloor (\bar{e}_0 \rfloor \tilde{d}\tilde{\theta}^3)\tilde{\theta}^3 \bigg) \Bigg] \\
= & -\frac{1}{2}\big[-\bar{e}_0 \rfloor (\bar{e}_0 \rfloor \tilde{d}\tilde{\theta}^0)\tilde{\theta}^1\big] = 0 \, .
\end{aligned}
$$

При дальнейшем вычислении будем использовать формулы (2.2.2). Откуда следует, что есть смысл рассматривать лишь суммационный индекс $d = 0$, так как при всех остальных внешние дифференциалы от базисных 1-форм (2.1.6) зануляются. Определив однозначно $d$ и используя метрику (2.1.5), можем утверждать, что последнее слагаемое не равно нулю только при суммационном индексе $= 1$ (**Приложение 2**)

Следовательно:

$$
\begin{aligned}
\tilde{\Gamma}_{12} &= H_x \tilde{\theta}^1, \quad \tilde{\Gamma}_{21} = -H_x \tilde{\theta}^1, \\
\tilde{\Gamma}_{13} &= H_y \tilde{\theta}^1, \quad \tilde{\Gamma}_{31} = -H_y \tilde{\theta}^1 \, .
\end{aligned} \qquad (2.2.4)
$$

Теперь получим вид 1-формы связности. Будем использовать известное соотношение:

$$
\tilde{\Gamma}^a{}_c = g^{ab} \tilde{\Gamma}_{bc} \qquad (2.2.5)
$$

$$\tilde{\Gamma}^a{}_1 = g^{ab}\tilde{\Gamma}_{b1} = g^{a0}\tilde{\Gamma}_{01} + g^{a1}\tilde{\Gamma}_{11} + g^{a2}\tilde{\Gamma}_{21} + g^{a3}\tilde{\Gamma}_{31}\,.$$

Для нахождения 1-формы связности переберем все возможные значения индекса $a$.

$$\tilde{\Gamma}^0{}_1 = g^{0b}\tilde{\Gamma}_{b1} = \underbrace{g^{00}}_{0}\underbrace{\tilde{\Gamma}_{01}}_{0} + g^{01}\underbrace{\tilde{\Gamma}_{11}}_{0} + \underbrace{g^{02}}_{0}\tilde{\Gamma}_{21} + \underbrace{g^{03}}_{0}\tilde{\Gamma}_{31} = 0\,,$$

$$\tilde{\Gamma}^1{}_1 = g^{10}\underbrace{\tilde{\Gamma}_{01}}_{0} + \underbrace{g^{11}}_{0}\underbrace{\tilde{\Gamma}_{11}}_{0} + \underbrace{g^{12}}_{0}\tilde{\Gamma}_{21} + \underbrace{g^{13}}_{0}\tilde{\Gamma}_{31} = 0\,,$$

$$\tilde{\Gamma}^2{}_1 = \underbrace{g^{20}}_{0}\tilde{\Gamma}_{01} + \underbrace{g^{21}}_{0}\tilde{\Gamma}_{11} + \underbrace{g^{22}}_{-1}\tilde{\Gamma}_{21} + \underbrace{g^{23}}_{0}\tilde{\Gamma}_{31} = -\tilde{\Gamma}_{21}\,,$$

$$\tilde{\Gamma}^3{}_1 = \underbrace{g^{30}}_{0}\tilde{\Gamma}_{01} + \underbrace{g^{31}}_{0}\tilde{\Gamma}_{11} + \underbrace{g^{32}}_{0}\tilde{\Gamma}_{21} + \underbrace{g^{33}}_{-1}\tilde{\Gamma}_{31} = -\tilde{\Gamma}_{31}\,,$$

$$\tilde{\Gamma}^a{}_2 = g^{ab}\tilde{\Gamma}_{b2} = g^{a0}\tilde{\Gamma}_{02} + g^{a1}\tilde{\Gamma}_{12} + g^{a2}\tilde{\Gamma}_{22} + g^{a3}\tilde{\Gamma}_{32}\,,$$

$$\tilde{\Gamma}^0{}_2 = \underbrace{g^{00}}_{0}\tilde{\Gamma}_{02} + \underbrace{g^{01}}_{1}\tilde{\Gamma}_{12} + \underbrace{g^{02}}_{0}\tilde{\Gamma}_{22} + \underbrace{g^{03}}_{0}\tilde{\Gamma}_{32} = \tilde{\Gamma}_{12}\,,$$

$$\tilde{\Gamma}^1{}_2 = \underbrace{g^{10}}_{1}\underbrace{\tilde{\Gamma}_{02}}_{0} + \underbrace{g^{11}}_{0}\tilde{\Gamma}_{12} + \underbrace{g^{12}}_{0}\tilde{\Gamma}_{22} + \underbrace{g^{13}}_{0}\tilde{\Gamma}_{32} = 0\,,$$

$$\tilde{\Gamma}^2{}_2 = \underbrace{g^{20}}_{0}\tilde{\Gamma}_{02} + \underbrace{g^{21}}_{0}\tilde{\Gamma}_{12} + \underbrace{g^{22}}_{0}\tilde{\Gamma}_{22} + \underbrace{g^{23}}_{0}\tilde{\Gamma}_{32} = 0\,,$$

$$\tilde{\Gamma}^3{}_2 = \underbrace{g^{30}}_{0}\tilde{\Gamma}_{02} + \underbrace{g^{31}}_{0}\tilde{\Gamma}_{12} + \underbrace{g^{32}}_{0}\tilde{\Gamma}_{22} + \underbrace{g^{33}}_{-1}\underbrace{\tilde{\Gamma}_{32}}_{0} = 0\,,$$

$$\tilde{\Gamma}^a{}_3 = g^{ab}\tilde{\Gamma}_{b3} = g^{a0}\tilde{\Gamma}_{03} + g^{01}\tilde{\Gamma}_{13} + g^{a2}\tilde{\Gamma}_{23} + g^{a3}\tilde{\Gamma}_{33}\,,$$

$$\tilde{\Gamma}^0{}_3 = \underbrace{g^{00}}_{0}\tilde{\Gamma}_{03} + \underbrace{g^{01}}_{1}\tilde{\Gamma}_{13} + \underbrace{g^{02}}_{0}\tilde{\Gamma}_{23} + \underbrace{g^{03}}_{0}\tilde{\Gamma}_{33} = \tilde{\Gamma}_{13}\,,$$

$$\tilde{\Gamma}^1{}_3 = \underbrace{g^{01}}_{1}\underbrace{\tilde{\Gamma}_{03}}_{0} + \underbrace{g^{11}}_{0}\tilde{\Gamma}_{13} + \underbrace{g^{12}}_{0}\tilde{\Gamma}_{23} + \underbrace{g^{13}}_{0}\tilde{\Gamma}_{33} = 0\,,$$

$$\tilde{\Gamma}^2{}_3 = \underbrace{g^{20}}_{0}\tilde{\Gamma}_{03} + \underbrace{g^{21}}_{0}\tilde{\Gamma}_{13} + \underbrace{g^{22}}_{-1}\underbrace{\tilde{\Gamma}_{23}}_{0} + \underbrace{g^{23}}_{0}\tilde{\Gamma}_{33} = 0\,,$$

$$\tilde{\Gamma}^3{}_3 = \underbrace{g^{30}}_{0}\tilde{\Gamma}_{03} + \underbrace{g^{31}}_{0}\tilde{\Gamma}_{13} + \underbrace{g^{32}}_{0}\tilde{\Gamma}_{23} + \underbrace{g^{33}}_{-1}\underbrace{\tilde{\Gamma}_{33}}_{0} = 0\,,$$

Таким образом, получим:

$$\begin{aligned}
\tilde{\Gamma}^0{}_2 &= H_x \tilde{d}u \\
\tilde{\Gamma}^0{}_3 &= H_y \tilde{d}u \\
\tilde{\Gamma}^2{}_1 &= H_x \tilde{d}u \\
\tilde{\Gamma}^3{}_1 &= H_y \tilde{d}u
\end{aligned}$$

(2.2.6)

Все остальные 1-формы связности равны нулю.

Далее получим явный вид коэффициентов связности. Напомним, что 1-формы связности связаны с коэффициентами связности следующим образом:

$$\tilde{\Gamma}^\alpha{}_\beta = \Gamma^\alpha{}_{\beta\gamma} \tilde{\theta}^\gamma \ ,$$

То есть $\Gamma^\alpha{}_{\beta\gamma}$ являются компонентами 1-формы связности $\tilde{\Gamma}^\alpha{}_\beta$.

Тогда раскроем суммационный индекс:

$\tilde{\Gamma}^0{}_2 = \Gamma^0{}_{2\gamma}\tilde{\theta}^\gamma = \Gamma^0{}_{20}\tilde{\theta}^0 + \Gamma^0{}_{21}\tilde{\theta}^1 + \Gamma^0{}_{22}\tilde{\theta}^2 + \Gamma^0{}_{23}\tilde{\theta}^3.$

С другой стороны, согласно (2.2.6) имеем: $\tilde{\Gamma}^0{}_2 = H_x \tilde{d}u.$

Приравниваем коэффициенты при $\tilde{\theta}^1$:

$\boldsymbol{\Gamma^0{}_{21} = H_x.}$

$\tilde{\Gamma}^0{}_3 = \Gamma^0{}_{3\gamma}\tilde{\theta}^\gamma = \Gamma^0{}_{30}\tilde{\theta}^0 + \Gamma^0{}_{31}\tilde{\theta}^1 + \Gamma^0{}_{32}\tilde{\theta}^2 + \Gamma^0{}_{33}\tilde{\theta}^3,$

С другой стороны по (2.2.6) имеем: $\tilde{\Gamma}^0{}_3 = H_y \tilde{d}u.$

Приравниваем коэффициенты при $\tilde{\theta}^1$:

$\boldsymbol{\Gamma^0{}_{31} = H_y.}$

$\tilde{\Gamma}^2{}_1 = \Gamma^2{}_{1\gamma}\tilde{\theta}^\gamma = \Gamma^2{}_{10}\tilde{\theta}^0 + \Gamma^2{}_{11}\tilde{\theta}^1 + \Gamma^2{}_{12}\tilde{\theta}^2 + \Gamma^2{}_{13}\tilde{\theta}^3.$

С другой стороны имеем: $\tilde{\Gamma}^2{}_1 = \tilde{d}u.$

Приравниваем коэффициенты при $\tilde{\theta}^1$:

$\boldsymbol{\Gamma^2{}_{11} = H_x.}$

$\tilde{\Gamma}^3{}_1 = \Gamma^3{}_{1\gamma}\tilde{\theta}^\gamma = \Gamma^3{}_{10}\tilde{\theta}^0 + \Gamma^3{}_{11}\tilde{\theta}^1 + \Gamma^3{}_{12}\tilde{\theta}^2 + \Gamma^3{}_{13}\tilde{\theta}^3.$

С другой стороны имеем: $\tilde{\Gamma}^3{}_1 = H_y \tilde{d}u$.

Приравниваем коэффициенты при $\tilde{\theta}^1$:

$$\Gamma^3{}_{11} = H_y.$$

Можем выписать коэффициенты связности в виде матриц следующим образом:

$$\Gamma^0{}_{\alpha\beta} = \begin{pmatrix} 0 & 0 & 0 & 0 \\ 0 & 0 & 0 & 0 \\ 0 & H_x & 0 & 0 \\ 0 & H_y & 0 & 0 \end{pmatrix} \qquad \Gamma^1{}_{\alpha\beta} = \begin{pmatrix} 0 & 0 & 0 & 0 \\ 0 & 0 & 0 & 0 \\ 0 & 0 & 0 & 0 \\ 0 & 0 & 0 & 0 \end{pmatrix}$$

$$\Gamma^2{}_{\alpha\beta} = \begin{pmatrix} 0 & 0 & 0 & 0 \\ 0 & H_x & 0 & 0 \\ 0 & 0 & 0 & 0 \\ 0 & 0 & 0 & 0 \end{pmatrix} \qquad \Gamma^3{}_{\alpha\beta} = \begin{pmatrix} 0 & 0 & 0 & 0 \\ 0 & H_y & 0 & 0 \\ 0 & 0 & 0 & 0 \\ 0 & 0 & 0 & 0 \end{pmatrix}$$

### §2.3 2-форма кривизны, тензор Римана и тензор Риччи для волновой метрики на многообразиях без кручения

Для вычисления 2-формы кривизны используется второе структурное уравнение Картана, рассмотренное в §10 первой главы. Это уравнение имеет вид:

$$\tilde{\mathcal{R}}^\alpha{}_\beta = \tilde{d}\tilde{\Gamma}^\alpha{}_\beta + \tilde{\Gamma}^\alpha{}_\gamma \wedge \tilde{\Gamma}^\gamma{}_\beta\ .$$

Данное структурное уравнение Картана можно расписать в следующем виде:

$$\tilde{\mathcal{R}}^\alpha{}_\beta = \tilde{d}\tilde{\Gamma}^\alpha{}_\beta + \tilde{\Gamma}^\alpha{}_0 \wedge \tilde{\Gamma}^0{}_\beta + \tilde{\Gamma}^\alpha{}_1 \wedge \tilde{\Gamma}^1{}_\beta + \tilde{\Gamma}^\alpha{}_2 \wedge \tilde{\Gamma}^2{}_\beta + \tilde{\Gamma}^\alpha{}_3 \wedge \tilde{\Gamma}^3{}_\beta\ .$$

Отличные от нуля выражения для 2-формы кривизны примут вид:

$$\boxed{\begin{array}{l} \tilde{\mathcal{R}}^0{}_2 = H_{xx}\tilde{\theta}^2 \wedge \tilde{\theta}^1 + H_{xy}\tilde{\theta}^3 \wedge \tilde{\theta}^1 \\ \tilde{\mathcal{R}}^0{}_3 = H_{yy}\tilde{\theta}^3 \wedge \tilde{\theta}^1 + H_{yx}\tilde{\theta}^2 \wedge \tilde{\theta}^1 \\ \tilde{\mathcal{R}}^2{}_1 = H_{xx}\tilde{\theta}^2 \wedge \tilde{\theta}^1 + H_{xy}\tilde{\theta}^3 \wedge \tilde{\theta}^1 \\ \tilde{\mathcal{R}}^3{}_1 = H_{yy}\tilde{\theta}^3 \wedge \tilde{\theta}^1 + H_{yx}\tilde{\theta}^2 \wedge \tilde{\theta}^1 \end{array}}$$

(2.3.1)

Вывод изложен в **Приложении 3**.

Далее, на основе разложения 2-формы по базису, рассчитаем компоненты тензора Римана:

$$\tilde{\mathcal{R}}^a{}_b = \frac{1}{2} R^a{}_{bcd} \tilde{\theta}^c \wedge \tilde{\theta}^d \overline{\overline{c<d}} \ R^a{}_{bcd} \tilde{\theta}^c \wedge \tilde{\theta}^d \ .$$

Тогда:

$$\tilde{\mathcal{R}}^0{}_2 \overline{\overline{c<d}} \ R^0{}_{2cd} \tilde{\theta}^c \wedge \tilde{\theta}^d$$
$$= R^0{}_{201} \tilde{\theta}^0 \wedge \tilde{\theta}^1 + R^2{}_{202} \tilde{\theta}^0 \wedge \tilde{\theta}^2 + R^0{}_{203} \tilde{\theta}^0 \wedge \tilde{\theta}^3 + R^0{}_{212} \tilde{\theta}^1 \wedge \tilde{\theta}^2$$
$$+ R^0{}_{213} \tilde{\theta}^1 \wedge \tilde{\theta}^3 + R^0{}_{223} \tilde{\theta}^2 \wedge \tilde{\theta}^3 \ .$$

С другой стороны, согласно (2.3.1):

$$\tilde{\mathcal{R}}^0{}_2 = H_{xx} \tilde{\theta}^2 \wedge \tilde{\theta}^1 + H_{xy} \tilde{\theta}^3 \wedge \tilde{\theta}^1 \ .$$

Приравниваем коэффициенты при $\tilde{\theta}^2 \wedge \tilde{\theta}^1$ и $\tilde{\theta}^3 \wedge \tilde{\theta}^1$, остальные нули. С учетом антисимметрии по двум первым и двум последним индексам, получаем:

$$-R^0{}_{212} = R^0{}_{221} = H_{xx} \ , \quad -R^0{}_{213} = R^0{}_{231} = H_{xy} \ .$$

Аналогично:

$$\tilde{\mathcal{R}}^0{}_3 \overline{\overline{c<d}} \ R^0{}_{3cd} \tilde{\theta}^c \wedge \tilde{\theta}^d$$
$$= R^0{}_{301} \tilde{\theta}^0 \wedge \tilde{\theta}^1 + R^0{}_{302} \tilde{\theta}^0 \wedge \tilde{\theta}^2 + R^0{}_{303} \tilde{\theta}^0 \wedge \tilde{\theta}^3 + R^0{}_{312} \tilde{\theta}^1 \wedge \tilde{\theta}^2$$
$$+ R^0{}_{313} \tilde{\theta}^1 \wedge \tilde{\theta}^3 + R^0{}_{323} \tilde{\theta}^2 \wedge \tilde{\theta}^2 \ .$$

С другой стороны:

$$\tilde{\mathcal{R}}^0{}_3 = H_{yx} \tilde{\theta}^2 \wedge \tilde{\theta}^1 + H_{yy} \tilde{\theta}^3 \wedge \tilde{\theta}^1 \ .$$

Приравниваем последовательно коэффициенты при $\tilde{\theta}^2 \wedge \tilde{\theta}^1$ и $\tilde{\theta}^3 \wedge \tilde{\theta}^1$, остальные нули. С учетом антисимметрии по последним индексам, получаем отличные от нуля коэффициента тензора Римана:

$$-R^0{}_{312} = R^0{}_{321} = H_{yx} \ , \quad -R^0{}_{313} = R^0{}_{331} = H_{yy} \ .$$

$$\tilde{\mathcal{R}}^2{}_1 \overline{\overline{c<d}} \ R^2{}_{1cd} \tilde{\theta}^c \wedge \tilde{\theta}^d$$
$$= R^2{}_{101} \tilde{\theta}^0 \wedge \tilde{\theta}^1 + R^2{}_{102} \tilde{\theta}^0 \wedge \tilde{\theta}^2 + R^2{}_{103} \tilde{\theta}^0 \wedge \tilde{\theta}^3 + R^2{}_{112} \tilde{\theta}^1 \wedge \tilde{\theta}^2$$
$$+ R^2{}_{113} \tilde{\theta}^1 \wedge \tilde{\theta}^3 + R^2{}_{123} \tilde{\theta}^2 \wedge \tilde{\theta}^3 \ .$$

С другой стороны:

$$\tilde{\mathcal{R}}^2{}_1 = H_{xx} \tilde{\theta}^2 \wedge \tilde{\theta}^1 + H_{xy} \tilde{\theta}^3 \wedge \tilde{\theta}^1 \ .$$

Приравниваем коэффициенты при $\tilde{\theta}^2\wedge\tilde{\theta}^1$ и $\tilde{\theta}^3\wedge\tilde{\theta}^1$, остальные нули:

$$-R^2{}_{112} = R^2{}_{121} = H_{xx}\ ,\quad R^2{}_{131} = -R^2{}_{113} = H_{xy}\ .$$

Далее:

$$\tilde{\mathcal{R}}^3{}_1 \overline{\overline{c<d}}\ R^3{}_{1cd}\tilde{\theta}^c\wedge\tilde{\theta}^d$$
$$= R^3{}_{101}\tilde{\theta}^0\wedge\tilde{\theta}^1 + R^3{}_{102}\tilde{\theta}^0\wedge\tilde{\theta}^2 + R^3{}_{103}\tilde{\theta}^0\wedge\tilde{\theta}^3 + R^3{}_{112}\tilde{\theta}^1\wedge\tilde{\theta}^2$$
$$+ R^3{}_{113}\tilde{\theta}^1\wedge\tilde{\theta}^3 + R^3{}_{123}\tilde{\theta}^2\wedge\tilde{\theta}^3\ .$$

С другой стороны:

$$\tilde{\mathcal{R}}^3{}_1 = H_{yx}\tilde{\theta}^2\wedge\tilde{\theta}^1 + H_{yy}\tilde{\theta}^3\wedge\tilde{\theta}^1\ .$$

Приравниваем коэффициенты при $\tilde{\theta}^2\wedge\tilde{\theta}^1$ и $\tilde{\theta}^3\wedge\tilde{\theta}^1$ соответственно, остальные нули. Получаем компоненты тензора Римана:

$$R^3{}_{121} = -R^3{}_{112} = H_{yx}\ ,\quad R^3{}_{131} = -R^3{}_{113} = H_{yy}\ .$$

Таким образом, получим следующие отличные от нуля компоненты тензора Римана:

$$\boxed{\begin{array}{l} R^0{}_{221} = -R^0{}_{212} = H_{xx} \\ R^0{}_{321} = -R^0{}_{312} = H_{yx} \\ R^0{}_{331} = -R^0{}_{313} = H_{yy} \\ R^0{}_{231} = -R^0{}_{213} = H_{xy} \\ R^2{}_{121} = -R^2{}_{112} = H_{xx} \\ R^2{}_{131} = -R^2{}_{113} = H_{xy} \\ R^3{}_{121} = -R^3{}_{112} = H_{yx} \\ R^3{}_{131} = -R^3{}_{113} = H_{yy} \end{array}}$$

(2.3.2)

Напомним вид уравнения Эйнштейна в пустом пространстве:

$$R_{ab} - \frac{1}{2}g_{ab}R = 0\ . \tag{2.3.3}$$

Правая часть уравнения равна нулю, так как тензор энергии-импульса материи в пустоте отсутствует. Как говорилось в первом параграфе данной главы, уравнение (2.3.3) в пустом пространстве имеет вид:

$$R_{ab} = 0\ . \tag{2.3.4}$$

Получим вид этого уравнения для метрики (2.1.4), используя вычисленные компоненты тензора Римана (2.3.2). Для вычисления тензора Риччи используем его определение:

$$R_{ab} = R^c{}_{acb} = R^c{}_{adb}\delta^d_c \ , \qquad (2.3.5)$$

Где $R^c{}_{adb}$ – компоненты тензора Риччи, вычисленные выше, а $\delta^d_c$ - символ Кронекера.

$$\begin{aligned}R_{00} = {}&R^0{}_{000}\delta^0_0 + R^0{}_{010}\delta^0_1 + R^0{}_{020}\delta^0_2 + R^0{}_{030}\delta^0_3 + R^1{}_{000}\delta^1_0 + R^1{}_{010}\delta^1_1 \\ &+ R^1{}_{020}\delta^1_2 + R^1{}_{030}\delta^1_3 + R^2{}_{000}\delta^2_0 + R^2{}_{010}\delta^2_1 + R^2{}_{020}\delta^2_2 + R^2{}_{030}\delta^2_3 \\ &+ R^3{}_{000}\delta^3_0 + R^3{}_{010}\delta^3_1 + R^3{}_{020}\delta^3_2 + R^3{}_{030}\delta^3_3 \ .\end{aligned}$$

Используя определение символа Кронекера и (2.3.2) получим следующий вид последнего равенства:

$$R_{00} = \underbrace{R^0{}_{000}}_{0} + \underbrace{R^1{}_{010}}_{0} + \underbrace{R^2{}_{020}}_{0} + \underbrace{R^3{}_{030}}_{0} = 0 \ .$$

Аналогично вычислим и все остальные компонента тензора Риччи. (**Приложение 4.**)

Итак, только одна компонента тензора Риччи отлична от нуля, а именно:

$$R_{11} = H_{xx} + H_{yy} \ . \qquad (2.3.6)$$

Теперь проверим, что скалярная кривизна $R$ для рассмотренной волновой метрики в многообразии без кручения равна нулю.

При доказательстве будем использовать определение $R$:

$$R = R_{ab}g^{ab} = R_{11}\underbrace{g^{11}}_{0} = 0 \ . \qquad (2.3.7)$$

Итак, в результате проведенных вычислений для волновой метрики (2.1.5) на многообразиях с нулевым кручением получили следующее уравнение:

$$H_{xx} + H_{yy} = 0 \ . \qquad (2.3.8)$$

или

$$\frac{\partial H}{\partial x^2} + \frac{\partial H}{\partial y^2} = 0 \ . \tag{2.3.9}$$

Уравнение (2.3.9) – это дифференциальное уравнение эллиптического типа, называемое уравнением Лапласа. Это уравнение (2.3.9) является искомым условием, при котором волновая метрика (2.1.5) является решением уравнения Эйнштейна в пустоте.

## Заключение

В работе было найдено условие, при котором волновая метрика риманова пространства является решением уравнения Эйнштейна в пустоте. Для этого сначала были изучены геометрические структуры на дифференцируемом многообразии: связность, кривизна и кручение связности. При этом используется такой аналитический аппарат дифференциальной геометрии, как исчисление внешних дифференциальных форм.

На основам дифференциальной геометрии был сделан реферативный обзор по прилагаемой литературе, все использованные в работе формулы были вычислены и изложены в единой системе обозначений. На основе внешнего дифференциального исчисления были выведены структурные уравнения Картана. Первое структурное уравнение, описывающее такую геометрическую структуру многообразия, как кручение связности, было выведено на основе понятия канонической векторнозначной 1-формы и свойства обобщенного внешнего дифференциала. При получении второго структурного уравнения Картана для 2-формы кривизны используется формула для обобщенного внешнего дифференциала от векторнозначной базисной 0-формы. На основе первого структурного уравнения Картана выводится формула, позволяющая определить 1-форму связности риманова пространства (при равном нулю кручении).

Далее, используя первое и второе структурные уравнения Картана для многообразия без кручения наделенного волновой метрикой, в работе были вычислены: 1-формы связности и коэффициенты связности, а также 2-формы кривизны, тензор Римана, тензор Риччи и скаляр кривизны.

На основе вычисленных геометрических характеристик многообразия получен вид уравнения Эйнштейна в пустоте. Это уравнение является дифференциальным уравнением эллиптического типа, называемое уравнением Лапласа и представляющее собой сумму вторых частных производных от функции волновой метрики по соответствующим координатам. Полученное уравнение является искомым условием, при котором рассматриваемая волновая метрика является решением уравнения Эйнштейна в пустоте.

## Литература

1. В.Т. Базылев. Геометрия дифференцируемых многообразий. – М.: Высш. шк., 1989. – 224 с.
2. Б.Шутц. Геометрические методы в математической физике. – М.: Платон, 1995. – 304 с.
3. В.И Арнольд. Математические методы классической механики. – М.: Едиториал УРСС, 2003. – 416 с.
4. В.Н. Пономарев, А.О. Барвинский, Ю.Н. Обухов. Геометродинамические методы и калибровочный подход к теории гравитационных взаимодействий. – М.: Энергоатомиздат, 1985. – 168 с.
5. А. П. Норден. Пространства аффинной связности. – М.: Наука, 1976. – 432 с.
6. Ч. Мизнер, К. Торн, Дж. Уилер. Гравитация Т.1-3. – М.: Едиториал УРСС, 1995. – 1512 с.
7. П.К. Рашевский. Риманова геометрия и тензорный анализ. – М.: Едиториал УРСС. 2006. – 664 с.
8. А. Бессе. Многообразия Эйнштейна. В 2-х томах. Т. 1-2. – М.: Мир, 1990. – 704 с.
9. В.Д. Захаров. Гравитационные волны в теории тяготения Эйнштейна. – М.: Наука, 1972.
10. С. Хелгасон. Дифференциальная геометрия, группы Ли и симметрические пространства. – М.: Факториал, 2005. – 605 с.
11. W. Adamowicz. Plane waves in gauge theories of gravitation.// Gen. Rel. Grav. – 1980. – V.12. – N9. – pp. 677-691.
12. Ehlers J., Kundt W. Exact solutions of the gravitational field equations.// Gravitation / ed L. Witten, John Wiley. – New York. – 1962. – P. 49-101.
13. Hehl F.W., McCrea J.D., Mielke E.W., Neeman Y. Metric affine gauge theory of gravity: field equations, Noether identities, world spinors and breaking of dilaton invariance. // Phys. Rep. – 1995. – V.258. – P. 1-171.
14. Курапова Е.И. Волновое решение уравнений Эйнштейна.

## Приложение 1

Рассмотрим многообразия, где метрика и связность согласованы, и получим выражение для коэффициентов связности через метрический тензор.

$$\nabla_\alpha g_{\beta\gamma} = \partial_\alpha g_{\beta\gamma} - \Gamma^\sigma_{\beta\alpha} g_{\sigma\gamma} - \Gamma^\sigma_{\gamma\alpha} g_{\beta\sigma} = 0$$

$$\nabla_\beta g_{\gamma\alpha} = \partial_\beta g_{\gamma\alpha} - \Gamma^\sigma_{\gamma\beta} g_{\sigma\alpha} - \Gamma^\sigma_{\alpha\beta} g_{\gamma\sigma} = 0$$

$$\nabla_\gamma g_{\alpha\beta} = \partial_\gamma g_{\alpha\beta} - \Gamma^\sigma_{\alpha\gamma} g_{\sigma\beta} - \Gamma^\sigma_{\beta\gamma} g_{\alpha\sigma} = 0$$

$$\left(g_{\alpha\beta} = \breve{g}(\bar{e}_\alpha, \bar{e}_\beta) = \bar{e}_\alpha \cdot \bar{e}_\beta = \bar{e}_\beta \cdot \bar{e}_\alpha = \breve{g}(\bar{e}_\beta, \bar{e}_\alpha) = g_{\beta\alpha}\right)$$

$$\nabla_\alpha g_{\beta\gamma} + \nabla_\beta g_{\gamma\alpha} - \nabla_\gamma g_{\alpha\beta} = \partial_\alpha g_{\beta\gamma} + \partial_\beta g_{\gamma\alpha} - \partial_\gamma g_{\alpha\beta} - \underbrace{\left(\Gamma^\sigma_{\beta\alpha} + \Gamma^\sigma_{\alpha\beta}\right)}_{2\Gamma^\sigma_{(\alpha\beta)}} g_{\sigma\gamma} +$$

$$+ \underbrace{\left(\Gamma^\sigma_{\alpha\gamma} - \Gamma^\sigma_{\gamma\alpha}\right)}_{2\Gamma^\sigma_{[\alpha\gamma]}} g_{\sigma\beta} + \underbrace{\left(\Gamma^\sigma_{\beta\gamma} - \Gamma^\sigma_{\gamma\beta}\right)}_{2\Gamma^\sigma_{[\beta\gamma]}} g_{\alpha\sigma} = 0$$

Многообразия, на которых антисимметричная часть связности равна нулю, называются **многообразиями Римана ($V_4$)**. В $V_4$ метрический тензор и связность согласованы, и связность симметрична (т.е. $\Gamma^\sigma_{\beta\gamma} = \Gamma^\sigma_{\gamma\beta}$).

Пространства, на которых антисимметричная часть связности отлична от нуля, то есть, **пространства с кручением ($U_4$)**, впервые были рассмотрены французским математиком Э. *Картаном*. На этих пространствах при параллельном переносе вектора сам контур не замыкается, и компоненты тензора кручения имеют вид: $2\Gamma^\sigma_{[\beta\gamma]} = T^\sigma_{\gamma\beta}$.

$$\partial_\alpha g_{\beta\gamma} + \partial_\beta g_{\gamma\alpha} - \partial_\gamma g_{\alpha\beta} - 2\Gamma^\sigma_{(\alpha\beta)} g_{\sigma\gamma} - 2\Gamma^\sigma_{[\alpha\beta]} g_{\sigma\gamma} + 2\Gamma^\sigma_{[\alpha\beta]} g_{\sigma\gamma} +$$

$$+ 2\Gamma^\sigma_{[\alpha\gamma]} g_{\sigma\beta} + 2\Gamma^\sigma_{[\beta\gamma]} g_{\alpha\sigma} = 0$$

$$\partial_\alpha g_{\beta\gamma} + \partial_\beta g_{\gamma\alpha} - \partial_\gamma g_{\alpha\beta} - 2\Gamma^\sigma_{\alpha\beta} g_{\sigma\gamma} - T^\sigma_{\alpha\beta} g_{\sigma\gamma} + T^\sigma_{\gamma\alpha} g_{\sigma\beta} - T^\sigma_{\beta\gamma} g_{\alpha\sigma} = 0$$

$$\Gamma^\sigma_{\alpha\beta} g_{\sigma\gamma} = \frac{1}{2}\left(\partial_\alpha g_{\beta\gamma} + \partial_\beta g_{\gamma\alpha} - \partial_\gamma g_{\alpha\beta} - T_{\alpha\beta\gamma} + T_{\beta\gamma\alpha} - T_{\gamma\alpha\beta}\right)$$

$$\Gamma^\sigma_{\alpha\beta} \underbrace{g_{\sigma\gamma} g^{\gamma\lambda}}_{\delta^\lambda_\sigma} = \frac{1}{2} g^{\gamma\lambda}\left(\partial_\alpha g_{\gamma\beta} + \partial_\beta g_{\alpha\gamma} - \partial_\gamma g_{\alpha\beta} + T_{\alpha\gamma\beta} - T_{\beta\alpha\gamma} + T_{\gamma\beta\alpha}\right)$$

введем обозначение: $\partial_\eta g_{\tau\varphi} = g_{\tau\varphi,\eta}$

$$\Gamma^\sigma_{\alpha\beta} = \frac{1}{2} g^{\lambda\gamma}\left(g_{\gamma\beta,\alpha} + g_{\alpha\gamma,\beta} - g_{\beta\alpha,\gamma}\right) + \frac{1}{2}\left(T_\alpha{}^\lambda{}_\beta + T^\lambda{}_{\beta\alpha} - T_{\beta\alpha}{}^\lambda\right)$$

$$\Gamma^\lambda_{\mu\nu} = \frac{1}{2} g^{\sigma\lambda}\left(g_{\sigma\nu,\mu} - g_{\nu\mu,\sigma} + g_{\mu\sigma,\nu}\right) + \frac{1}{2}\left(T^\lambda{}_{\nu\mu} - T_{\nu\mu}{}^\lambda + T_\mu{}^\lambda{}_\nu\right)$$

$$\Gamma^{\lambda}_{\mu\nu} = \frac{1}{2}g^{\sigma\lambda}\big(\partial_{\mu}g_{\nu\sigma} + \partial_{\nu}g_{\sigma\mu} - \partial_{\sigma}g_{\mu\nu} - T_{\mu\nu\sigma} + T_{\nu\sigma\mu} - T_{\sigma\mu\nu}\big)$$

$$\Gamma^{\lambda}_{\mu\nu} = \frac{1}{2}g^{\lambda\sigma}\Delta^{\alpha\beta\gamma}_{\sigma\nu\mu}\big(g_{\alpha\beta,\gamma} + T_{\alpha\beta\gamma}\big)$$

$\Delta^{\alpha\beta\gamma}_{\sigma\nu\mu} = \delta^{\alpha}_{\sigma}\delta^{\beta}_{\nu}\delta^{\gamma}_{\mu} - \delta^{\alpha}_{\nu}\delta^{\beta}_{\mu}\delta^{\gamma}_{\sigma} + \delta^{\alpha}_{\mu}\delta^{\beta}_{\sigma}\delta^{\gamma}_{\nu}$ – скобка Схоутэна.

**Приложение 2**

$$\tilde{\Gamma}_{01} = -\frac{1}{2}\left[g_{1d}(\bar{e}_0\rfloor\tilde{d}\tilde{\theta}^d) - g_{0d}(\bar{e}_1\rfloor\tilde{d}\tilde{\theta}^d) - g_{cd}\bar{e}_0\rfloor(\bar{e}_1\rfloor\tilde{d}\tilde{\theta}^d)\tilde{\theta}^c\right]$$

$$= -\frac{1}{2}\left[g_{1d}(\bar{e}_0\rfloor\tilde{d}\tilde{\theta}^d) - g_{0d}(\bar{e}_1\rfloor\tilde{d}\tilde{\theta}^d) - g_{0d}\bar{e}_0\rfloor(\bar{e}_1\rfloor\tilde{d}\tilde{\theta}^d)\tilde{\theta}^0\right.$$

$$\left. - g_{1d}\bar{e}_0\rfloor(\bar{e}_1\rfloor\tilde{d}\tilde{\theta}^d)\tilde{\theta}^1 - g_{2d}\bar{e}_0\rfloor(\bar{e}_1\rfloor\tilde{d}\tilde{\theta}^d)\tilde{\theta}^2 - g_{3d}\bar{e}_0\rfloor(\bar{e}_1\rfloor\tilde{d}\tilde{\theta}^d)\tilde{\theta}^3\right]$$

$$= -\frac{1}{2}\left[\underbrace{g_{10}}_{1}(\bar{e}_0\rfloor\tilde{d}\tilde{\theta}^0) - \underbrace{g_{00}}_{0}(\bar{e}_1\rfloor\tilde{d}\tilde{\theta}^0) - \underbrace{g_{00}}_{0}\bar{e}_0\rfloor(\bar{e}_1\rfloor\tilde{d}\tilde{\theta}^0)\tilde{\theta}^0\right.$$

$$\left. - \underbrace{g_{10}}_{1}\bar{e}_0\rfloor(\bar{e}_1\rfloor\tilde{d}\tilde{\theta}^0)\tilde{\theta}^1 - \underbrace{g_{20}}_{0}\bar{e}_0\rfloor(\bar{e}_1\rfloor\tilde{d}\tilde{\theta}^0)\tilde{\theta}^2 - \underbrace{g_{30}}_{0}\bar{e}_0\rfloor(\bar{e}_1\rfloor\tilde{d}\tilde{\theta}^0)\tilde{\theta}^3\right]$$

$$= -\frac{1}{2}\left[\bar{e}_0\rfloor\tilde{d}\tilde{\theta}^0 - \bar{e}_0\rfloor(\bar{e}_1\rfloor\tilde{d}\tilde{\theta}^0)\tilde{\theta}^1\right] = 0.$$

Аналогично вычисляя, найдем значения для других 1-форм связности:

$$\tilde{\Gamma}_{02} = -\frac{1}{2}\left[g_{2d}(\bar{e}_0\rfloor\tilde{d}\tilde{\theta}^d) - g_{0d}(\bar{e}_2\rfloor\tilde{d}\tilde{\theta}^d) - g_{cd}\bar{e}_0\rfloor(\bar{e}_2\rfloor\tilde{d}\tilde{\theta}^d)\tilde{\theta}^c\right]$$

$$= -\frac{1}{2}\left[g_{2d}(\bar{e}_0\rfloor\tilde{d}\tilde{\theta}^d) - g_{0d}(\bar{e}_2\rfloor\tilde{d}\tilde{\theta}^d) - g_{0d}\bar{e}_0\rfloor(\bar{e}_2\rfloor\tilde{d}\tilde{\theta}^d)\tilde{\theta}^0\right.$$

$$\left. - g_{1d}\bar{e}_0\rfloor(\bar{e}_2\rfloor\tilde{d}\tilde{\theta}^d)\tilde{\theta}^1 - g_{2d}\bar{e}_0\rfloor(\bar{e}_2\rfloor\tilde{d}\tilde{\theta}^d)\tilde{\theta}^2 - g_{3d}\bar{e}_0\rfloor(\bar{e}_2\rfloor\tilde{d}\tilde{\theta}^d)\tilde{\theta}^3\right]$$

$$= -\frac{1}{2}\left[\underbrace{g_{20}}_{0}(\bar{e}_0\rfloor\tilde{d}\tilde{\theta}^0) - \underbrace{g_{00}}_{0}(\bar{e}_2\rfloor\tilde{d}\tilde{\theta}^0) - \underbrace{g_{00}}_{0}\bar{e}_0\rfloor(\bar{e}_2\rfloor\tilde{d}\tilde{\theta}^0)\tilde{\theta}^0\right.$$

$$\left. - \underbrace{g_{10}}_{1}\bar{e}_0\rfloor(\bar{e}_2\rfloor\tilde{d}\tilde{\theta}^0)\tilde{\theta}^1\right] = -\frac{1}{2}\left[-\bar{e}_0\rfloor(\bar{e}_2\rfloor\tilde{d}\tilde{\theta}^0)\tilde{\theta}^1\right] = 0,$$

$$\tilde{\Gamma}_{03} = -\frac{1}{2}\left[g_{3d}(\bar{e}_0\rfloor\tilde{d}\tilde{\theta}^d) - g_{0d}(\bar{e}_3\rfloor\tilde{d}\tilde{\theta}^d) - g_{cd}\bar{e}_0\rfloor(\bar{e}_3\rfloor\tilde{d}\tilde{\theta}^d)\tilde{\theta}^c\right]$$

$$= -\frac{1}{2}\left[g_{3d}(\bar{e}_0\rfloor\tilde{d}\tilde{\theta}^d) - g_{0d}(\bar{e}_3\rfloor\tilde{d}\tilde{\theta}^d) - g_{0d}\bar{e}_0\rfloor(\bar{e}_3\rfloor\tilde{d}\tilde{\theta}^d)\tilde{\theta}^0\right.$$

$$\left. - g_{1d}\bar{e}_0\rfloor(\bar{e}_3\rfloor\tilde{d}\tilde{\theta}^d)\tilde{\theta}^1 - g_{2d}\bar{e}_0\rfloor(\bar{e}_3\rfloor\tilde{d}\tilde{\theta}^d)\tilde{\theta}^2 - g_{3d}\bar{e}_0\rfloor(\bar{e}_3\rfloor\tilde{d}\tilde{\theta}^d)\tilde{\theta}^3\right]$$

$$= -\frac{1}{2}\left[\underbrace{g_{30}}_{0}(\bar{e}_0\rfloor\tilde{d}\tilde{\theta}^0) - \underbrace{g_{00}}_{0}(\bar{e}_3\rfloor\tilde{d}\tilde{\theta}^0) - \underbrace{g_{10}}_{1}\bar{e}_0\rfloor(\bar{e}_3\rfloor\tilde{d}\tilde{\theta}^0)\tilde{\theta}^1\right]$$

$$= -\frac{1}{2}\left[-\bar{e}_0\rfloor(\bar{e}_3\rfloor\tilde{d}\tilde{\theta}^0)\tilde{\theta}^1\right] = 0,$$

$$\tilde{\Gamma}_{10} = -\frac{1}{2}[g_{0d}(\bar{e}_1 \rfloor \tilde{d}\tilde{\theta}^d) - g_{1d}(\bar{e}_0 \rfloor \tilde{d}\tilde{\theta}^d) - g_{cd}\bar{e}_1 \rfloor (\bar{e}_0 \rfloor \tilde{d}\tilde{\theta}^d)\tilde{\theta}^c]$$

$$= -\frac{1}{2}[g_{0d}(\bar{e}_1 \rfloor \tilde{d}\tilde{\theta}^d) - g_{1d}(\bar{e}_0 \rfloor \tilde{d}\tilde{\theta}^d) - g_{0d}\bar{e}_1 \rfloor (\bar{e}_0 \rfloor \tilde{d}\tilde{\theta}^d)\tilde{\theta}^0$$

$$- g_{1d}\bar{e}_1 \rfloor (\bar{e}_0 \rfloor \tilde{d}\tilde{\theta}^d)\tilde{\theta}^1 - g_{2d}\bar{e}_1 \rfloor (\bar{e}_0 \rfloor \tilde{d}\tilde{\theta}^d)\tilde{\theta}^2 - g_{3d}\bar{e}_1 \rfloor (\bar{e}_0 \rfloor \tilde{d}\tilde{\theta}^d)\tilde{\theta}^3]$$

$$= -\frac{1}{2}\left[\underbrace{g_{00}}_{0}(\bar{e}_1 \rfloor \tilde{d}\tilde{\theta}^0) - \underbrace{g_{10}}_{1}(\bar{e}_0 \rfloor \tilde{d}\tilde{\theta}^0) - \underbrace{g_{10}}_{1}\bar{e}_1 \rfloor (\bar{e}_0 \rfloor \tilde{d}\tilde{\theta}^0)\tilde{\theta}^1\right]$$

$$= -\frac{1}{2}[-\bar{e}_0 \rfloor \tilde{d}\tilde{\theta}^0 - \bar{e}_1 \rfloor (\bar{e}_0 \rfloor \tilde{d}\tilde{\theta}^0)\tilde{\theta}^1] = 0,$$

$$\tilde{\Gamma}_{11} = -\frac{1}{2}[g_{1d}(\bar{e}_1 \rfloor \tilde{d}\tilde{\theta}^d) - g_{1d}(\bar{e}_1 \rfloor \tilde{d}\tilde{\theta}^d) - g_{cd}\bar{e}_1 \rfloor (\bar{e}_1 \rfloor \tilde{d}\tilde{\theta}^d)\tilde{\theta}^0].$$

В последней формуле третье слагаемое зануляется, в силу (2.2.2).

Тогда формула принимает вид:

$$\tilde{\Gamma}_{11} = -\frac{1}{2}[g_{1d}(\bar{e}_1 \rfloor \tilde{d}\tilde{\theta}^d) - g_{1d}(\bar{e}_1 \rfloor \tilde{d}\tilde{\theta}^d)] = -\frac{1}{2}\left[\underbrace{g_{10}}_{1}(\bar{e}_1 \rfloor \tilde{d}\tilde{\theta}^0) - \underbrace{g_{10}}_{1}(\bar{e}_1 \rfloor \tilde{d}\tilde{\theta}^0)\right]$$

$$= -\frac{1}{2}[\bar{e}_1 \rfloor \tilde{d}\tilde{\theta}^0 - \bar{e}_1 \rfloor \tilde{d}\tilde{\theta}^0] = 0,$$

$$\tilde{\Gamma}_{12} = -\frac{1}{2}[g_{2d}(\bar{e}_1 \rfloor \tilde{d}\tilde{\theta}^d) - g_{1d}(\bar{e}_2 \rfloor \tilde{d}\tilde{\theta}^d) - g_{cd}\bar{e}_1 \rfloor (\bar{e}_2 \rfloor \tilde{d}\tilde{\theta}^d)\tilde{\theta}^c]$$

$$= -\frac{1}{2}[g_{2d}(\bar{e}_1 \rfloor \tilde{d}\tilde{\theta}^d) - g_{1d}(\bar{e}_2 \rfloor \tilde{d}\tilde{\theta}^d) - g_{0d}\bar{e}_1 \rfloor (\bar{e}_2 \rfloor \tilde{d}\tilde{\theta}^d)\tilde{\theta}^0$$

$$- g_{1d}\bar{e}_1 \rfloor (\bar{e}_2 \rfloor \tilde{d}\tilde{\theta}^d)\tilde{\theta}^1 - g_{2d}\bar{e}_1 \rfloor (\bar{e}_2 \rfloor \tilde{d}\tilde{\theta}^d)\tilde{\theta}^2 - g_{3d}\bar{e}_1 \rfloor (\bar{e}_2 \rfloor \tilde{d}\tilde{\theta}^1)\tilde{\theta}^3]$$

$$= -\frac{1}{2}\left[\underbrace{g_{20}}_{0}(\bar{e}_1 \rfloor \tilde{d}\tilde{\theta}^0) - \underbrace{g_{10}}_{1}(\bar{e}_2 \rfloor \tilde{d}\tilde{\theta}^0) - \underbrace{g_{10}}_{1}\bar{e}_1 \rfloor (\bar{e}_2 \rfloor \tilde{d}\tilde{\theta}^0)\tilde{\theta}^1\right]$$

$$= -\frac{1}{2}[-\bar{e}_2 \rfloor \tilde{d}\tilde{\theta}^0 - \bar{e}_1 \rfloor (\bar{e}_2 \rfloor \tilde{d}\tilde{\theta}^0)\tilde{\theta}^1],$$

$$\tilde{\Gamma}_{13} = -\frac{1}{2}[g_{3d}(\bar{e}_1 \rfloor \tilde{d}\tilde{\theta}^d) - g_{1d}(\bar{e}_3 \rfloor \tilde{d}\tilde{\theta}^d) - g_{cd}\bar{e}_1 \rfloor (\bar{e}_3 \rfloor \tilde{d}\tilde{\theta}^d)\tilde{\theta}^c]$$

$$= -\frac{1}{2}[g_{3d}(\bar{e}_1 \rfloor \tilde{d}\tilde{\theta}^d) - g_{1d}(\bar{e}_3 \rfloor \tilde{d}\tilde{\theta}^d) - g_{0d}\bar{e}_1 \rfloor (\bar{e}_3 \rfloor \tilde{d}\tilde{\theta}^d)\tilde{\theta}^0$$

$$- g_{1d}\bar{e}_1 \rfloor (\bar{e}_3 \rfloor \tilde{d}\tilde{\theta}^d)\tilde{\theta}^1 - g_{2d}\bar{e}_1 \rfloor (\bar{e}_3 \rfloor \tilde{d}\tilde{\theta}^d)\tilde{\theta}^2 - g_{3d}\bar{e}_1 \rfloor (\bar{e}_3 \rfloor \tilde{d}\tilde{\theta}^d)\tilde{\theta}^3]$$

$$= -\frac{1}{2}\left[\underbrace{g_{30}}_{0}(\bar{e}_1 \rfloor \tilde{d}\tilde{\theta}^0) - \underbrace{g_{10}}_{1}(\bar{e}_3 \rfloor \tilde{d}\tilde{\theta}^0) - \underbrace{g_{10}}_{1}\bar{e}_1 \rfloor (\bar{e}_3 \rfloor \tilde{d}\tilde{\theta}^0)\tilde{\theta}^1\right]$$

$$= -\frac{1}{2}[-\bar{e}_3 \rfloor \tilde{d}\tilde{\theta}^0 - \bar{e}_1 \rfloor (\bar{e}_3 \rfloor \tilde{d}\tilde{\theta}^0)\tilde{\theta}^1],$$

$$\tilde{\Gamma}_{20} = -\frac{1}{2}\big[g_{0d}(\bar{e}_2\rfloor \tilde{d}\tilde{\theta}^d) - g_{2d}(\bar{e}_0\rfloor \tilde{d}\tilde{\theta}^d) - g_{cd}\bar{e}_2\rfloor(\bar{e}_0\rfloor \tilde{d}\tilde{\theta}^d)\tilde{\theta}^c\big]$$

$$= -\frac{1}{2}\big[g_{0d}(\bar{e}_2\rfloor \tilde{d}\tilde{\theta}^d) - g_{2d}(\bar{e}_0\rfloor \tilde{d}\tilde{\theta}^d) - g_{0d}\bar{e}_2\rfloor(\bar{e}_0\rfloor \tilde{d}\tilde{\theta}^d)\tilde{\theta}^0$$

$$- g_{1d}\bar{e}_2\rfloor(\bar{e}_0\rfloor \tilde{d}\tilde{\theta}^d)\tilde{\theta}^1 - g_{2d}(\bar{e}_0\rfloor \tilde{d}\tilde{\theta}^d)\tilde{\theta}^2 - g_{3d}(\bar{e}_0\rfloor \tilde{d}\tilde{\theta}^d)\tilde{\theta}^3\big]$$

$$= -\frac{1}{2}\bigg[\underbrace{g_{00}}_{0}(\bar{e}_2\rfloor \tilde{d}\tilde{\theta}^0) - \underbrace{g_{20}}_{0}(\bar{e}_0\rfloor \tilde{d}\tilde{\theta}^0) - \underbrace{g_{10}}_{1}\bar{e}_2\rfloor(\bar{e}_0\rfloor \tilde{d}\tilde{\theta}^0)\tilde{\theta}^1\bigg]$$

$$= -\frac{1}{2}\big[-\bar{e}_2\rfloor(\bar{e}_0\rfloor \tilde{d}\tilde{\theta}^0)\tilde{\theta}^1\big] = 0,$$

$$\tilde{\Gamma}_{21} = -\frac{1}{2}\big[g_{1d}(\bar{e}_2\rfloor \tilde{d}\tilde{\theta}^d) - g_{2d}(\bar{e}_1\rfloor \tilde{d}\tilde{\theta}^d) - g_{cd}\bar{e}_2\rfloor(\bar{e}_1\rfloor \tilde{d}\tilde{\theta}^d)\tilde{\theta}^c\big]$$

$$= -\frac{1}{2}\big[g_{1d}(\bar{e}_2\rfloor \tilde{d}\tilde{\theta}^d) - g_{2d}(\bar{e}_1\rfloor \tilde{d}\tilde{\theta}^d) - g_{0d}\bar{e}_2\rfloor(\bar{e}_1\rfloor \tilde{d}\tilde{\theta}^d)\tilde{\theta}^0$$

$$- g_{1d}\bar{e}_2\rfloor(\bar{e}_1\rfloor \tilde{d}\tilde{\theta}^d)\tilde{\theta}^1 - g_{2d}\bar{e}_2\rfloor(\bar{e}_1\rfloor \tilde{d}\tilde{\theta}^d)\tilde{\theta}^2 - g_{3d}\bar{e}_2\rfloor(\bar{e}_1\rfloor \tilde{d}\tilde{\theta}^d)\tilde{\theta}^3\big]$$

$$= -\frac{1}{2}\bigg[\underbrace{g_{10}}_{1}(\bar{e}_2\rfloor \tilde{d}\tilde{\theta}^0) - \underbrace{g_{20}}_{0}(\bar{e}_1\rfloor \tilde{d}\tilde{\theta}^0) - \underbrace{g_{10}}_{1}\bar{e}_2\rfloor(\bar{e}_1\rfloor \tilde{d}\tilde{\theta}^0)\tilde{\theta}^1\bigg]$$

$$= -\frac{1}{2}\big[\bar{e}_2\rfloor \tilde{d}\tilde{\theta}^0 - \bar{e}_2\rfloor(\bar{e}_1\rfloor \tilde{d}\tilde{\theta}^0)\tilde{\theta}^1\big].$$

Рассматривая следующую формулу, сразу заметим, что последнее слагаемое равно нулю (используем (2.2.2)):

$$\tilde{\Gamma}_{22} = -\frac{1}{2}\big[g_{2d}(\bar{e}_2\rfloor \tilde{d}\tilde{\theta}^d) - g_{2d}(\bar{e}_2\rfloor \tilde{d}\tilde{\theta}^d) - g_{cd}\bar{e}_2\rfloor(\bar{e}_2\rfloor \tilde{d}\tilde{\theta}^d)\tilde{\theta}^c\big]$$

$$= -\frac{1}{2}\big[g_{2d}(\bar{e}_2\rfloor \tilde{d}\tilde{\theta}^d) - g_{2d}(\bar{e}_2\rfloor \tilde{d}\tilde{\theta}^d)\big]$$

$$= -\frac{1}{2}\bigg[\underbrace{g_{20}}_{0}(\bar{e}_2\rfloor \tilde{d}\tilde{\theta}^0) - \underbrace{g_{20}}_{0}(\bar{e}_2\rfloor \tilde{d}\tilde{\theta}^0)\bigg] = 0,$$

$$\tilde{\Gamma}_{23} = -\frac{1}{2}\big[g_{3d}(\bar{e}_2\rfloor \tilde{d}\tilde{\theta}^d) - g_{2d}(\bar{e}_3\rfloor \tilde{d}\tilde{\theta}^d) - g_{cd}\bar{e}_2\rfloor(\bar{e}_3\rfloor \tilde{d}\tilde{\theta}^d)\tilde{\theta}^c\big]$$

$$= -\frac{1}{2}\big[g_{3d}(\bar{e}_2\rfloor \tilde{d}\tilde{\theta}^d) - g_{2d}(\bar{e}_3\rfloor \tilde{d}\tilde{\theta}^d) - g_{0d}\bar{e}_2\rfloor(\bar{e}_3\rfloor \tilde{d}\tilde{\theta}^d)\tilde{\theta}^0$$

$$- g_{1d}\bar{e}_2\rfloor(\bar{e}_3\rfloor \tilde{d}\tilde{\theta}^d)\tilde{\theta}^1 - g_{2d}\bar{e}_2\rfloor(\bar{e}_3\rfloor \tilde{d}\tilde{\theta}^d)\tilde{\theta}^2 - g_{3d}\bar{e}_2\rfloor(\bar{e}_3\rfloor \tilde{d}\tilde{\theta}^d)\tilde{\theta}^3\big].$$

Зная, что только $\tilde{d}\tilde{\theta}^0 \neq 0$, получим:

$$\tilde{\Gamma}_{23} = -\frac{1}{2}\bigg[\underbrace{g_{03}}_{0}(\bar{e}_2\rfloor \tilde{d}\tilde{\theta}^0) - \underbrace{g_{20}}_{0}(\bar{e}_3\rfloor \tilde{d}\tilde{\theta}^0) - \underbrace{g_{10}}_{1}\bar{e}_2\rfloor(\bar{e}_3\rfloor \tilde{d}\tilde{\theta}^0)\tilde{\theta}^1\bigg]$$

$$= -\frac{1}{2}\big[-\bar{e}_2\rfloor(\bar{e}_3\rfloor \tilde{d}\tilde{\theta}^0)\tilde{\theta}^1\big] = 0,$$

$$\tilde{\Gamma}_{30} = -\frac{1}{2}\big[g_{0d}(\bar{e}_3\rfloor\tilde{d}\tilde{\theta}^d) - g_{3d}(\bar{e}_0\rfloor\tilde{d}\tilde{\theta}^d) - g_{cd}\bar{e}_3\rfloor(\bar{e}_0\rfloor\tilde{d}\tilde{\theta}^d)\tilde{\theta}^c\big]$$

$$= -\frac{1}{2}\big[g_{0d}(\bar{e}_3\rfloor\tilde{d}\tilde{\theta}^d) - g_{3d}(\bar{e}_0\rfloor\tilde{d}\tilde{\theta}^d) - g_{0d}\bar{e}_3\rfloor(\bar{e}_0\rfloor\tilde{d}\tilde{\theta}^d)\tilde{\theta}^0$$

$$- g_{1d}\bar{e}_3\rfloor(\bar{e}_0\rfloor\tilde{d}\tilde{\theta}^d)\tilde{\theta}^1 - g_{2d}\bar{e}_3\rfloor(\bar{e}_0\rfloor\tilde{d}\tilde{\theta}^d)\tilde{\theta}^2 - g_{3d}\bar{e}_3\rfloor(\bar{e}_0\rfloor\tilde{d}\tilde{\theta}^d)\tilde{\theta}^3\big]$$

$$= -\frac{1}{2}\Big[\underbrace{g_{00}}_{0}(\bar{e}_3\rfloor\tilde{d}\tilde{\theta}^0) - \underbrace{g_{30}}_{0}(\bar{e}_0\rfloor\tilde{d}\tilde{\theta}^0) - \underbrace{g_{10}}_{1}\bar{e}_3\rfloor(\bar{e}_0\rfloor\tilde{d}\tilde{\theta}^0)\tilde{\theta}^1\Big]$$

$$= -\frac{1}{2}\big[-\bar{e}_3\rfloor(\bar{e}_0\rfloor\tilde{d}\tilde{\theta}^0)\tilde{\theta}^1\big] = 0\,,$$

$$\tilde{\Gamma}_{31} = -\frac{1}{2}\big[g_{1d}(\bar{e}_3\rfloor\tilde{d}\tilde{\theta}^d) - g_{3d}(\bar{e}_1\rfloor\tilde{d}\tilde{\theta}^d) - g_{cd}\bar{e}_3\rfloor(\bar{e}_1\rfloor\tilde{d}\tilde{\theta}^d)\tilde{\theta}^c\big]$$

$$= -\frac{1}{2}\big[g_{1d}(\bar{e}_3\rfloor\tilde{d}\tilde{\theta}^d) - g_{3d}(\bar{e}_1\rfloor\tilde{d}\tilde{\theta}^d) - g_{0d}\bar{e}_3\rfloor(\bar{e}_1\rfloor\tilde{d}\tilde{\theta}^d)\tilde{\theta}^0$$

$$- g_{1d}\bar{e}_3\rfloor(\bar{e}_1\rfloor\tilde{d}\tilde{\theta}^d)\tilde{\theta}^1 - g_{2d}\bar{e}_3\rfloor(\bar{e}_1\rfloor\tilde{d}\tilde{\theta}^d)\tilde{\theta}^2 - g_{3d}\bar{e}_3\rfloor(\bar{e}_1\rfloor\tilde{d}\tilde{\theta}^d)\tilde{\theta}^3\big]$$

$$= -\frac{1}{2}\Big[\underbrace{g_{10}}_{1}(\bar{e}_3\rfloor\tilde{d}\tilde{\theta}^0) - \underbrace{g_{30}}_{0}(\bar{e}_1\rfloor\tilde{d}\tilde{\theta}^0) - \underbrace{g_{10}}_{1}\bar{e}_3\rfloor(\bar{e}_1\rfloor\tilde{d}\tilde{\theta}^0)\tilde{\theta}^1\Big]$$

$$= -\frac{1}{2}\big[\bar{e}_3\rfloor\tilde{d}\tilde{\theta}^0 - \bar{e}_3\rfloor(\bar{e}_1\rfloor\tilde{d}\tilde{\theta}^0)\tilde{\theta}^1\big]\,,$$

$$\tilde{\Gamma}_{32} = -\frac{1}{2}\big[g_{2d}(\bar{e}_3\rfloor\tilde{d}\tilde{\theta}^d) - g_{3d}(\bar{e}_2\rfloor\tilde{d}\tilde{\theta}^d) - g_{cd}\bar{e}_3\rfloor(\bar{e}_2\rfloor\tilde{d}\tilde{\theta}^d)\tilde{\theta}^c\big]$$

$$= -\frac{1}{2}\big[g_{2d}(\bar{e}_3\rfloor\tilde{d}\tilde{\theta}^d) - g_{3d}(\bar{e}_2\rfloor\tilde{d}\tilde{\theta}^d) - g_{0d}\bar{e}_3\rfloor(\bar{e}_2\rfloor\tilde{d}\tilde{\theta}^d)\tilde{\theta}^0$$

$$- g_{1d}\bar{e}_3\rfloor(\bar{e}_2\rfloor\tilde{d}\tilde{\theta}^d)\tilde{\theta}^1 - g_{2d}\bar{e}_3\rfloor(\bar{e}_2\rfloor\tilde{d}\tilde{\theta}^d)\tilde{\theta}^2 - g_{3d}\bar{e}_3\rfloor(\bar{e}_2\rfloor\tilde{d}\tilde{\theta}^d)\tilde{\theta}^3\big]$$

$$= -\frac{1}{2}\Big[\underbrace{g_{20}}_{0}(\bar{e}_3\rfloor\tilde{d}\tilde{\theta}^0) - \underbrace{g_{30}}_{0}(\bar{e}_2\rfloor\tilde{d}\tilde{\theta}^0) - \underbrace{g_{10}}_{1}\bar{e}_3\rfloor(\bar{e}_2\rfloor\tilde{d}\tilde{\theta}^0)\tilde{\theta}^1\Big]$$

$$= -\frac{1}{2}\big[-\bar{e}_3\rfloor(\bar{e}_2\rfloor\tilde{d}\tilde{\theta}^0)\tilde{\theta}^1\big] = 0\,.$$

Третье слагаемое в следующей сумме зануляется (в силу (2.2.2)):

$$\tilde{\Gamma}_{33} = -\frac{1}{2}\big[g_{3d}(\bar{e}_3\rfloor\tilde{d}\tilde{\theta}^d) - g_{3d}(\bar{e}_3\rfloor\tilde{d}\tilde{\theta}^d) - g_{cd}\bar{e}_3\rfloor(\bar{e}_3\rfloor\tilde{d}\tilde{\theta}^d)\tilde{\theta}^c\big]$$

$$= -\frac{1}{2}\big[g_{3d}(\bar{e}_3\rfloor\tilde{d}\tilde{\theta}^d) - g_{3d}(\bar{e}_3\rfloor\tilde{d}\tilde{\theta}^d) - g_{0d}\bar{e}_3\rfloor(\bar{e}_3\rfloor\tilde{d}\tilde{\theta}^d)\tilde{\theta}^0$$

$$- g_{1d}\bar{e}_3\rfloor(\bar{e}_3\rfloor\tilde{d}\tilde{\theta}^d)\tilde{\theta}^1 - g_{2d}\bar{e}_3\rfloor(\bar{e}_3\rfloor\tilde{d}\tilde{\theta}^d)\tilde{\theta}^2 - g_{3d}\bar{e}_3\rfloor(\bar{e}_3\rfloor\tilde{d}\tilde{\theta}^d)\tilde{\theta}^3\big]$$

$$= -\frac{1}{2}\Big[\underbrace{g_{30}}_{0}(\bar{e}_3\rfloor\tilde{d}\tilde{\theta}^0) - \underbrace{g_{30}}_{0}(\bar{e}_3\rfloor\tilde{d}\tilde{\theta}^0)\Big] = 0\,.$$

Итак, мы получили следующие ненулевые связности:

$$\tilde{\Gamma}_{12} = -\frac{1}{2}\left[-\bar{e}_2\rfloor\tilde{d}\tilde{\theta}^0 - \bar{e}_1\rfloor(\bar{e}_2\rfloor\tilde{d}\tilde{\theta}^0)\tilde{\theta}^1\right],$$

$$\tilde{\Gamma}_{13} = -\frac{1}{2}\left[-\bar{e}_3\rfloor\tilde{d}\tilde{\theta}^0 - \bar{e}_1\rfloor(\bar{e}_3\rfloor\tilde{d}\tilde{\theta}^0)\tilde{\theta}^1\right],$$

$$\tilde{\Gamma}_{21} = -\frac{1}{2}\left[\bar{e}_2\rfloor\tilde{d}\tilde{\theta}^0 - \bar{e}_2\rfloor(\bar{e}_0\rfloor\tilde{d}\tilde{\theta}^0)\tilde{\theta}^1\right],$$

$$\tilde{\Gamma}_{31} = -\frac{1}{2}\left[\bar{e}_3\rfloor\tilde{d}\tilde{\theta}^0 - \bar{e}_3\rfloor(\bar{e}_1\rfloor\tilde{d}\tilde{\theta}^0)\tilde{\theta}^1\right].$$

Все остальные связности при вычислениях занулятся.

Используя формулы (2.2.2), найдем значения этих связностей:

$$\tilde{\Gamma}_{12} = -\frac{1}{2}\left[-\bar{e}_2\rfloor\tilde{d}\tilde{\theta}^0 - \bar{e}_1\rfloor(\bar{e}_2\rfloor\tilde{d}\tilde{\theta}^0)\tilde{\theta}^1\right]$$
$$= \frac{1}{2}\left[\bar{e}_2\rfloor(-H_x\tilde{\theta}^1\wedge\tilde{\theta}^2 - H_y\tilde{\theta}^1\wedge\tilde{\theta}^3)\right.$$
$$\left. + \bar{e}_1\rfloor\left(\bar{e}_2\rfloor(-H_x\tilde{\theta}^1\wedge\tilde{\theta}^2 - H_y\tilde{\theta}^1\wedge\tilde{\theta}^3)\right)\tilde{\theta}^1\right]$$
$$= \frac{1}{2}\left[-H_x 2\delta_2^{[1}\tilde{\theta}^{2]} - H_y 2\delta_2^{[1}\tilde{\theta}^{3]} + \bar{e}_1\rfloor\left(H_x 2\delta_2^{[1}\tilde{\theta}^{2]} - H_y 2\delta_2^{[1}\tilde{\theta}^{3]}\right)\tilde{\theta}^1\right]$$
$$= \frac{1}{2}\left[-H_x 2\frac{1}{2}(\delta_2^1\tilde{\theta}^2 - \delta_2^2\tilde{\theta}^1) - H_y 2\frac{1}{2}(\delta_2^1\tilde{\theta}^3 - \delta_2^3\tilde{\theta}^1)\right.$$
$$\left. + \bar{e}_1\rfloor(H_x\tilde{\theta}^1)\tilde{\theta}^1\right] = \frac{1}{2}\left[H_x\tilde{\theta}^1 + H_x\bar{e}_1\rfloor\tilde{\theta}^1\tilde{\theta}^1\right] = H_x\tilde{\theta}^1.$$

Аналогично вычислим следующие связности:

$$\tilde{\Gamma}_{13} = -\frac{1}{2}\left[-\bar{e}_3\rfloor\tilde{d}\tilde{\theta}^0 - \bar{e}_1\rfloor(\bar{e}_3\rfloor\tilde{d}\tilde{\theta}^0)\tilde{\theta}^1\right]$$
$$= -\frac{1}{2}\left[-\bar{e}_3\rfloor(-H_x\tilde{\theta}^1\wedge\tilde{\theta}^2 - H_y\tilde{\theta}^1\wedge\tilde{\theta}^3)\right.$$
$$\left. - \bar{e}_1\rfloor\left(\bar{e}_3\rfloor(-H_x\tilde{\theta}^1\wedge\tilde{\theta}^2 - H_y\tilde{\theta}^1\wedge\tilde{\theta}^3)\right)\tilde{\theta}^1\right]$$
$$= -\frac{1}{2}\left[-H_x 2\delta_3^{[1}\tilde{\theta}^{2]} - H_y 2\delta_3^{[1}\tilde{\theta}^{3]}\right.$$
$$\left. - \bar{e}_1\rfloor\left(-H_x 2\delta_3^{[1}\tilde{\theta}^{2]} - H_y 2\delta_3^{[1}\tilde{\theta}^{3]}\right)\tilde{\theta}^1\right] = H_y\tilde{\theta}^1,$$

$$\tilde{\Gamma}_{21} = -\frac{1}{2}\left[\bar{e}_2\rfloor\tilde{d}\tilde{\theta}^0 - \bar{e}_2\rfloor(\bar{e}_0\rfloor\tilde{d}\tilde{\theta}^0)\tilde{\theta}^1\right]$$
$$= -\frac{1}{2}\left[\bar{e}_2\rfloor(-H_x\tilde{\theta}^1\wedge\tilde{\theta}^2 - H_y\tilde{\theta}^1\wedge\tilde{\theta}^3)\right.$$
$$\left. - \bar{e}_2\rfloor\left(\bar{e}_0\rfloor(-H_x\tilde{\theta}^1\wedge\tilde{\theta}^2 - H_y\tilde{\theta}^1\wedge\tilde{\theta}^3)\right)\tilde{\theta}^1\right] = -\frac{1}{2}\left[H_x\tilde{\theta}^1 + H_x\tilde{\theta}^1\right]$$
$$= -H_x\tilde{\theta}^1.$$

$$\tilde{\Gamma}_{31} = -\frac{1}{2}\big[\bar{e}_3 \rfloor \tilde{d}\tilde{\theta}^0 - \bar{e}_3\rfloor(\bar{e}_1\rfloor \tilde{d}\tilde{\theta}^0)\tilde{\theta}^1\big]$$
$$= -\frac{1}{2}\Big[\bar{e}_1\rfloor(-H_x\tilde{\theta}^1 \wedge \tilde{\theta}^2 - H_y\tilde{\theta}^1 \wedge \tilde{\theta}^3)$$
$$- \bar{e}_3\rfloor\big(\bar{e}_1\rfloor(-H_x\tilde{\theta}^1\wedge\tilde{\theta}^2 - H_y\tilde{\theta}^1\wedge\tilde{\theta}^3)\big)\tilde{\theta}^1\Big] = -H_y\tilde{\theta}^1 .$$

Следовательно:

$$\tilde{\Gamma}_{12} = H_x\tilde{\theta}^1, \quad \tilde{\Gamma}_{21} = -H_x\tilde{\theta}^1,$$
$$\tilde{\Gamma}_{13} = H_y\tilde{\theta}^1, \quad \tilde{\Gamma}_{31} = -H_y\tilde{\theta}^1 . \qquad (2.2.4)$$

## Приложение 3

Согласно (2.2.6):

$$\tilde{\mathcal{R}}^0{}_0 = \tilde{d}\underbrace{\tilde{\Gamma}^0{}_0}_{0} + \underbrace{\tilde{\Gamma}^0{}_0 \wedge \tilde{\Gamma}^0{}_0}_{0} + \tilde{\Gamma}^0{}_1 \wedge \tilde{\Gamma}^1{}_0 + \tilde{\Gamma}^0{}_2 \wedge \underbrace{\tilde{\Gamma}^2{}_0}_{0} + \tilde{\Gamma}^0{}_3 \wedge \underbrace{\tilde{\Gamma}^3{}_0}_{0} = 0 \ ,$$

$$\tilde{\mathcal{R}}^0{}_1 = \tilde{d}\underbrace{\tilde{\Gamma}^0{}_1}_{0} + \underbrace{\tilde{\Gamma}^0{}_0}_{0} \wedge \underbrace{\tilde{\Gamma}^0{}_1}_{0} + \underbrace{\tilde{\Gamma}^0{}_1}_{0} \wedge \underbrace{\tilde{\Gamma}^1{}_1}_{0} + \tilde{\Gamma}^0{}_2 \wedge \tilde{\Gamma}^2{}_1 + \tilde{\Gamma}^0{}_3 \wedge \tilde{\Gamma}^3{}_1 =$$
$$= \tilde{\Gamma}^0{}_2 \wedge \tilde{\Gamma}^2{}_1 + \tilde{\Gamma}^0{}_3 \wedge \tilde{\Gamma}^3{}_1 \ ,$$

$$\tilde{\mathcal{R}}^0{}_2 = \tilde{d}\tilde{\Gamma}^0{}_2 + \underbrace{\tilde{\Gamma}^0{}_0}_{0} \wedge \tilde{\Gamma}^0{}_2 + \underbrace{\tilde{\Gamma}^0{}_1}_{0} \wedge \underbrace{\tilde{\Gamma}^1{}_2}_{0} + \tilde{\Gamma}^0{}_2 \wedge \underbrace{\tilde{\Gamma}^2{}_2}_{0} + \tilde{\Gamma}^0{}_3 \wedge \underbrace{\tilde{\Gamma}^3{}_2}_{0} = \tilde{d}\tilde{\Gamma}^0{}_2 \ ,$$

$$\tilde{\mathcal{R}}^0{}_3 = \tilde{d}\tilde{\Gamma}^0{}_3 + \underbrace{\tilde{\Gamma}^0{}_0}_{0} \wedge \tilde{\Gamma}^0{}_3 + \underbrace{\tilde{\Gamma}^0{}_1}_{0} \wedge \underbrace{\tilde{\Gamma}^1{}_3}_{0} + \tilde{\Gamma}^0{}_2 \wedge \underbrace{\tilde{\Gamma}^2{}_3}_{0} + \tilde{\Gamma}^0{}_3 \wedge \underbrace{\tilde{\Gamma}^3{}_3}_{0} = \tilde{d}\tilde{\Gamma}^0{}_3 \ ,$$

$$\tilde{\mathcal{R}}^1{}_0 = \tilde{d}\underbrace{\tilde{\Gamma}^1{}_0}_{0} + \underbrace{\tilde{\Gamma}^1{}_0}_{0} \wedge \underbrace{\tilde{\Gamma}^0{}_0}_{0} + \underbrace{\tilde{\Gamma}^1{}_1}_{0} \wedge \underbrace{\tilde{\Gamma}^1{}_0}_{0} + \underbrace{\tilde{\Gamma}^1{}_2}_{0} \wedge \underbrace{\tilde{\Gamma}^2{}_0}_{0} + \underbrace{\tilde{\Gamma}^1{}_3}_{0} \wedge \underbrace{\tilde{\Gamma}^3{}_0}_{0} = 0 \ ,$$

$$\tilde{\mathcal{R}}^1{}_1 = \tilde{d}\underbrace{\tilde{\Gamma}^1{}_1}_{0} + \underbrace{\tilde{\Gamma}^1{}_0}_{0} \wedge \underbrace{\tilde{\Gamma}^0{}_1}_{0} + \underbrace{\tilde{\Gamma}^1{}_1 \wedge \tilde{\Gamma}^1{}_1}_{0} + \underbrace{\tilde{\Gamma}^1{}_2}_{0} \wedge \tilde{\Gamma}^2{}_1 + \underbrace{\tilde{\Gamma}^1{}_3}_{0} \wedge \tilde{\Gamma}^3{}_1 = 0 \ ,$$

$$\tilde{\mathcal{R}}^1{}_2 = \tilde{d}\underbrace{\tilde{\Gamma}^1{}_2}_{0} + \underbrace{\tilde{\Gamma}^1{}_0}_{0} \wedge \tilde{\Gamma}^0{}_2 + \underbrace{\tilde{\Gamma}^1{}_1}_{0} \wedge \underbrace{\tilde{\Gamma}^1{}_2}_{0} + \underbrace{\tilde{\Gamma}^1{}_2}_{0} \wedge \underbrace{\tilde{\Gamma}^2{}_2}_{0} + \underbrace{\tilde{\Gamma}^1{}_3}_{0} \wedge \underbrace{\tilde{\Gamma}^3{}_2}_{0} = 0 \ ,$$

$$\tilde{\mathcal{R}}^1{}_3 = \tilde{d}\underbrace{\tilde{\Gamma}^1{}_3}_{0} + \underbrace{\tilde{\Gamma}^1{}_0}_{0} \wedge \tilde{\Gamma}^0{}_3 + \underbrace{\tilde{\Gamma}^1{}_1}_{0} \wedge \underbrace{\tilde{\Gamma}^1{}_3}_{0} + \underbrace{\tilde{\Gamma}^1{}_2}_{0} \wedge \underbrace{\tilde{\Gamma}^2{}_3}_{0} + \underbrace{\tilde{\Gamma}^1{}_3}_{0} \wedge \underbrace{\tilde{\Gamma}^3{}_3}_{0} = 0 \ ,$$

$$\tilde{\mathcal{R}}^2{}_0 = \tilde{d}\underbrace{\tilde{\Gamma}^2{}_0}_{0} + \tilde{\Gamma}^2{}_0 \wedge \underbrace{\tilde{\Gamma}^0{}_0}_{0} + \tilde{\Gamma}^2{}_1 \wedge \underbrace{\tilde{\Gamma}^1{}_0}_{0} + \underbrace{\tilde{\Gamma}^2{}_2}_{0} \wedge \tilde{\Gamma}^2{}_0 + \underbrace{\tilde{\Gamma}^2{}_3}_{0} \wedge \tilde{\Gamma}^3{}_0 = 0 \ ,$$

$$\tilde{\mathcal{R}}^2{}_1 = \tilde{d}\tilde{\Gamma}^2{}_1 + \tilde{\Gamma}^2{}_0 \wedge \underbrace{\tilde{\Gamma}^0{}_1}_{0} + \tilde{\Gamma}^2{}_1 \wedge \underbrace{\tilde{\Gamma}^1{}_1}_{0} + \underbrace{\tilde{\Gamma}^2{}_2}_{0} \wedge \tilde{\Gamma}^2{}_1 + \underbrace{\tilde{\Gamma}^2{}_3}_{0} \wedge \underbrace{\tilde{\Gamma}^3{}_1}_{0} = \tilde{d}\tilde{\Gamma}^2{}_1 \ ,$$

$$\tilde{\mathcal{R}}^2{}_2 = \tilde{d}\underbrace{\tilde{\Gamma}^2{}_2}_{0} + \underbrace{\tilde{\Gamma}^2{}_0}_{0} \wedge \tilde{\Gamma}^0{}_2 + \tilde{\Gamma}^2{}_1 \wedge \underbrace{\tilde{\Gamma}^1{}_2}_{0} + \underbrace{\tilde{\Gamma}^2{}_2 \wedge \tilde{\Gamma}^2{}_2}_{0} + \underbrace{\tilde{\Gamma}^2{}_3}_{0} \wedge \underbrace{\tilde{\Gamma}^3{}_2}_{0} = 0 \ ,$$

$$\tilde{\mathcal{R}}^2{}_3 = \tilde{d}\underbrace{\tilde{\Gamma}^2{}_3}_{0} + \underbrace{\tilde{\Gamma}^2{}_0}_{0} \wedge \tilde{\Gamma}^0{}_3 + \tilde{\Gamma}^2{}_1 \wedge \underbrace{\tilde{\Gamma}^1{}_3}_{0} + \underbrace{\tilde{\Gamma}^2{}_2}_{0} \wedge \underbrace{\tilde{\Gamma}^2{}_3}_{0} + \underbrace{\tilde{\Gamma}^2{}_3}_{0} \wedge \underbrace{\tilde{\Gamma}^3{}_3}_{0} = 0 \ ,$$

$$\tilde{\mathcal{R}}^3{}_0 = \tilde{d}\underbrace{\tilde{\Gamma}^3{}_0}_{0} + \underbrace{\tilde{\Gamma}^3{}_0}_{0} \wedge \underbrace{\tilde{\Gamma}^0{}_0}_{0} + \tilde{\Gamma}^3{}_1 \wedge \underbrace{\tilde{\Gamma}^1{}_0}_{0} + \underbrace{\tilde{\Gamma}^3{}_2}_{0} \wedge \underbrace{\tilde{\Gamma}^2{}_0}_{0} + \underbrace{\tilde{\Gamma}^3{}_3}_{0} \wedge \underbrace{\tilde{\Gamma}^3{}_0}_{0} = 0 \ ,$$

$$\tilde{\mathcal{R}}^3{}_1 = \tilde{d}\tilde{\Gamma}^3{}_1 + \tilde{\Gamma}^3{}_0 \wedge \underbrace{\tilde{\Gamma}^0{}_1}_{0} + \tilde{\Gamma}^3{}_1 \wedge \underbrace{\tilde{\Gamma}^1{}_1}_{0} + \underbrace{\tilde{\Gamma}^3{}_2}_{0} \wedge \tilde{\Gamma}^2{}_1 + \underbrace{\tilde{\Gamma}^3{}_2}_{0} \wedge \tilde{\Gamma}^3{}_1 = \tilde{d}\tilde{\Gamma}^3{}_1 \ ,$$

$$\tilde{\mathcal{R}}^3{}_2 = \tilde{d}\underbrace{\tilde{\Gamma}^3{}_2}_{0} + \underbrace{\tilde{\Gamma}^3{}_0}_{0} \wedge \tilde{\Gamma}^0{}_2 + \tilde{\Gamma}^3{}_1 \wedge \underbrace{\tilde{\Gamma}^1{}_2}_{0} + \underbrace{\tilde{\Gamma}^3{}_2}_{0} \wedge \underbrace{\tilde{\Gamma}^2{}_2}_{0} + \underbrace{\tilde{\Gamma}^3{}_3}_{0} \wedge \underbrace{\tilde{\Gamma}^3{}_2}_{0} = 0 \ ,$$

$$\tilde{\mathcal{R}}^3{}_3 = \tilde{d}\underbrace{\tilde{\Gamma}^3{}_3}_{0} + \underbrace{\tilde{\Gamma}^3{}_0}_{0} \wedge \tilde{\Gamma}^0{}_3 + \tilde{\Gamma}^3{}_1 \wedge \underbrace{\tilde{\Gamma}^1{}_3}_{0} + \underbrace{\tilde{\Gamma}^3{}_2}_{0} \wedge \underbrace{\tilde{\Gamma}^2{}_3}_{0} + \underbrace{\tilde{\Gamma}^3{}_3 \wedge \tilde{\Gamma}^3{}_3}_{0} = 0 \ ,$$

Подставим (2.2.6) в полученные выражения для 2-формы кривизны:

$\tilde{\mathcal{R}}^0{}_1 = \tilde{\Gamma}^0{}_2 \wedge \tilde{\Gamma}^2{}_1 + \tilde{\Gamma}^0{}_3 \wedge \tilde{\Gamma}^3{}_1 = H_x \tilde{d}u \wedge H_x \tilde{d}u + H_y \tilde{d}u \wedge H_y \tilde{d}u = H_x{}^2 \tilde{d}u \wedge \tilde{d}u + H_y{}^2 \tilde{d}u \wedge \tilde{d}u = 0$ ,

слагаемые равны нулю по свойству оператора внешнего дифференцирования.

$$\begin{aligned}\tilde{\mathcal{R}}^0{}_2 &= \tilde{d}\tilde{\Gamma}^0{}_2 = \tilde{d}(H_x \tilde{d}u) = \partial_\alpha H_x \tilde{\theta}^\alpha \wedge \tilde{d}u \\ &= \partial_0 H_x \tilde{\theta}^0 \wedge \tilde{\theta}^1 + \partial_1 H_x \tilde{\theta}^1 \wedge \tilde{\theta}^1 + \partial_2 H_x \tilde{\theta}^2 \wedge \tilde{\theta}^1 + \partial_3 H_x \tilde{\theta}^3 \wedge \tilde{\theta}^1 \\ &= H_{xx} \tilde{\theta}^2 \wedge \tilde{\theta}^1 + H_{xy} \tilde{\theta}^3 \wedge \tilde{\theta}^1 \ , \end{aligned}$$

$\tilde{\mathcal{R}}^0{}_3 = \tilde{d}\tilde{\Gamma}^0{}_3 = \tilde{d}(H_y \tilde{d}u) = \partial_\alpha H_y \tilde{\theta}^\alpha \wedge \tilde{d}u = H_{yy} \tilde{\theta}^3 \wedge \tilde{\theta}^1 + H_{yx} \tilde{\theta}^2 \wedge \tilde{\theta}^1$ ,

$\tilde{\mathcal{R}}^2{}_1 = \tilde{d}\tilde{\Gamma}^2{}_1 = \tilde{d}(H_x \tilde{d}u) = H_{xx} \tilde{\theta}^2 \wedge \tilde{\theta}^1 + H_{xy} \tilde{\theta}^3 \wedge \tilde{\theta}^1$ ,

$\tilde{\mathcal{R}}^3{}_1 = \tilde{d}\tilde{\Gamma}^3{}_1 = \tilde{d}(H_y \tilde{d}u) = H_{yy} \tilde{\theta}^3 \wedge \tilde{\theta}^1 + H_{yx} \tilde{\theta}^2 \wedge \tilde{\theta}^1$ ,

где $H_{xx}, H_{yy}, H_{yx}, H_{xy}$ - соответствующие вторые производные.

Таким образом, получили следующие отличные от нуля выражения для 2-формы кривизны:

$$\boxed{\begin{aligned} \tilde{\mathcal{R}}^0{}_2 &= H_{xx} \tilde{\theta}^2 \wedge \tilde{\theta}^1 + H_{xy} \tilde{\theta}^3 \wedge \tilde{\theta}^1 \\ \tilde{\mathcal{R}}^0{}_3 &= H_{yy} \tilde{\theta}^3 \wedge \tilde{\theta}^1 + H_{yx} \tilde{\theta}^2 \wedge \tilde{\theta}^1 \\ \tilde{\mathcal{R}}^2{}_1 &= H_{xx} \tilde{\theta}^2 \wedge \tilde{\theta}^1 + H_{xy} \tilde{\theta}^3 \wedge \tilde{\theta}^1 \\ \tilde{\mathcal{R}}^3{}_1 &= H_{yy} \tilde{\theta}^3 \wedge \tilde{\theta}^1 + H_{yx} \tilde{\theta}^2 \wedge \tilde{\theta}^1 \end{aligned}}$$

(2.3.1)

## Приложение 4.

Так как при $c \neq d$ символ Кронекера $\delta_c^d = 0$, то для упрощения записи нулевые элементы выписывать не будем.

$$R_{01} = R^c{}_{0d1}\delta_c^d = \underbrace{R^0{}_{001}}_{0} + \underbrace{R^1{}_{011}}_{0} + \underbrace{R^2{}_{021}}_{0} + \underbrace{R^3{}_{031}}_{0} = 0\ ,$$

$$R_{02} = R^c{}_{0d2}\delta_c^d = \underbrace{R^0{}_{002}}_{0} + \underbrace{R^1{}_{012}}_{0} + \underbrace{R^2{}_{022}}_{0} + \underbrace{R^3{}_{032}}_{0} = 0\ ,$$

$$R_{03} = R^c{}_{0d3}\delta_c^d = \underbrace{R^0{}_{003}}_{0} + \underbrace{R^1{}_{013}}_{0} + \underbrace{R^2{}_{023}}_{0} + \underbrace{R^3{}_{033}}_{0} = 0\ ,$$

$$R_{10} = R^c{}_{1d0}\delta_c^d = \underbrace{R^0{}_{100}}_{0} + \underbrace{R^1{}_{110}}_{0} + \underbrace{R^2{}_{120}}_{0} + \underbrace{R^3{}_{130}}_{0} = 0\ ,$$

$$R_{11} = R^c{}_{1d1}\delta_c^d = \underbrace{R^0{}_{101}}_{0} + \underbrace{R^1{}_{111}}_{0} + R^2{}_{121} + R^3{}_{131} = H_{xx} + H_{yy}\ ,$$

$$R_{12} = R^c{}_{1d2}\delta_c^d = \underbrace{R^0{}_{102}}_{0} + \underbrace{R^1{}_{112}}_{0} + \underbrace{R^2{}_{122}}_{0} + \underbrace{R^3{}_{132}}_{0} = 0\ ,$$

$$R_{13} = R^c{}_{2d0}\delta_c^d = \underbrace{R^0{}_{201}}_{0} + \underbrace{R^1{}_{211}}_{0} + \underbrace{R^2{}_{221}}_{0} + \underbrace{R^3{}_{231}}_{0} = 0\ ,$$

$$R_{20} = R^c{}_{2d0}\delta_c^d = \underbrace{R^0{}_{200}}_{0} + \underbrace{R^1{}_{210}}_{0} + \underbrace{R^2{}_{220}}_{0} + \underbrace{R^3{}_{230}}_{0} = 0\ ,$$

$$R_{21} = R^c{}_{2d1}\delta_c^d = \underbrace{R^0{}_{201}}_{0} + \underbrace{R^1{}_{211}}_{0} + \underbrace{R^2{}_{221}}_{0} + \underbrace{R^3{}_{231}}_{0} = 0\ ,$$

$$R_{22} = R^c{}_{2d2}\delta_c^d = \underbrace{R^0{}_{202}}_{0} + \underbrace{R^1{}_{212}}_{0} + \underbrace{R^2{}_{222}}_{0} + \underbrace{R^3{}_{232}}_{0} = 0\ ,$$

$$R_{23} = R^c{}_{2d3}\delta_c^d = \underbrace{R^0{}_{203}}_{0} + \underbrace{R^1{}_{213}}_{0} + \underbrace{R^2{}_{223}}_{0} + \underbrace{R^3{}_{233}}_{0} = 0\ ,$$

$$R_{30} = R^c{}_{3d0}\delta_c^d = \underbrace{R^0{}_{300}}_{0} + \underbrace{R^1{}_{310}}_{0} + \underbrace{R^2{}_{320}}_{0} + \underbrace{R^3{}_{330}}_{0} = 0\ ,$$

$$R_{31} = R^c{}_{3d1}\delta_c^d = \underbrace{R^0{}_{301}}_{0} + \underbrace{R^1{}_{311}}_{0} + \underbrace{R^2{}_{321}}_{0} + \underbrace{R^3{}_{331}}_{0} = 0\ ,$$

$$R_{32} = R^c{}_{3d2}\delta_c^d = \underbrace{R^0{}_{302}}_{0} + \underbrace{R^1{}_{312}}_{0} + \underbrace{R^2{}_{322}}_{0} + \underbrace{R^3{}_{332}}_{0} = 0\ ,$$

$$R_{33} = R^c{}_{3d3}\delta_c^d = \underbrace{R^0{}_{303}}_{0} + \underbrace{R^1{}_{313}}_{0} + \underbrace{R^2{}_{323}}_{0} + \underbrace{R^3{}_{333}}_{0} = 0\ .$$